\documentclass[a4paper]{article}
\title{\bf Sobolev estimates for optimal transport maps \\
on Gaussian spaces}
\author{Shizan Fang$^a$\footnote{fang@u-bourgogne.fr}\quad Vincent Nolot$^a$\footnote{vincent.nolot@u-bourgogne.fr}
\vspace{3mm}\\
{\footnotesize $^a$I.M.B, BP 47870, Universit\'e de Bourgogne,
Dijon, France}
 }
\date{}
\usepackage{amssymb,amsmath,amsfonts,amsthm,color}

\def\R{\mathbb{R}}
\def\E{\mathbb{E}}
\def\P{\mathcal{P}}
\def\D{\mathbb{D}}

\def\Id{\hbox{\rm Id}}
\def\F{\mathcal{F}}
\def\L{\mathcal{L}}

\newcommand{\ra}{\rightarrow}

\newcommand{\ee}{\varepsilon}

\def\det{\textup{det}}
\def\Cylin{\textup{Cylin}}
\def\Ent{\textup{Ent}}
\def\div{\textup{div}}

\def\fin{\hfill$\square$}
\def\dis{\displaystyle}
\newtheorem{theorem}{Theorem}[section]
\newtheorem{lemma}[theorem]{Lemma}       
\newtheorem{corollary}[theorem]{Corollary}

\begin{document}

\maketitle
\makeatletter 
\renewcommand\theequation{\thesection.\arabic{equation}}
\@addtoreset{equation}{section}
\makeatother 

\begin{abstract}
We will study variations in Sobolev spaces of optimal transport maps with the standard Gaussian measure as the reference measure. 
Some dimension free inequalities will be obtained. As application, we construct solutions to Monge-Amp\`ere equations in finite dimension,
as well as on the Wiener space.
\end{abstract}
\vskip 2mm

 \vspace{2mm}


\vskip 10mm

\quad Let $e^{-V}dx$ and $e^{-W}dx$ be two probability measures on $\R^d$ having second moment, then there is a convex function $\Phi$ such that
$\nabla\Phi$ is the optimal transport map which pushes $e^{-V}dx$ to $e^{-W}dx$. If moreover (i) the functions $V$ and $W$ are smooth,
 bounded from below, (ii) the Hessian $\nabla^2V$ of $V$ is bounded from above and $\nabla W\geq K_1\,\Id$ with $K_1>0$, then 
$\Phi$ is smooth (see \cite{Caffarelli,Kolesnikov1}) and 
\begin{equation*}
\sup_{x\in\R^d}||\nabla^2\Phi(x)||_{HS}<+\infty,
\end{equation*}
where $||\cdot||_{HS}$ denotes the Hilbert-Schmidt norm. The above upper bound is dimension-dependent.
 In a recent work \cite{ Kolesnikov1}, A.V. Kolesnikov proved the inequality
\begin{equation}\label{0.1}
\int_{\R^d} |\nabla V|^2e^{-V}\,dx\geq K_1\,\int_{\R^d}||\nabla^2\Phi||_{HS}^2\,e^{-V}dx.
\end{equation}

Although the constant $K_1$ in \eqref{0.1} is of dimension free, but on infinite dimensional spaces, $\nabla^2\Phi$ usually is not of Hilbert-Schmidt class. Let 
$\nabla\Phi(x)=x+\nabla\varphi(x)$. A dimension free inequality for $||\nabla^2\varphi||_{HS}^2$ has been established in \cite{ Kolesnikov1}
 under the hypothesis 
 \begin{equation}\label{0.2}
 \nabla^2W\leq K_2\,\Id.
 \end{equation}
 
 Our work has been inspired from a series of works by A.V. Kolesnikov \cite{Kolesnikov1,Kolesnikov2,BogachevKolesnikov} 
 and a series of works by D. Feyel and A. S. \"Ust\"unel \cite{FeyelUstunel1,FeyelUstunel2,FeyelUstunel3}.
 The main contribution is to remove the condition \eqref{0.2}. Here is the result:
 
 \begin{theorem}\label{thA} Let $e^{-V}d\gamma$ and $e^{-W}d\gamma$ be two probability measures on $\R^d$, where $\gamma$ is the standard 
 Gaussian measure on $\R^d$. Suppose that $\nabla^2W\geq -c\,\Id$ with $c\in [0,1[$. Then
\begin{equation}\label{0.3}
\begin{split}
&\int_{\R^d}|\nabla V|^2 e^{-V}d\gamma-\int_{\R^d}|\nabla W|^2 e^{-W}d\gamma +{2\over 1-c}\int_{\R^d}||\nabla^2W||_{HS}^2\,e^{-W}d\gamma\\
&\geq 2\Ent_\gamma(e^{-V})-2\Ent_\gamma(e^{-W})+{1-c\over2}\int_{\R^d}||\nabla^2\varphi||_{HS}^2 \,e^{-V}d\gamma.
\end{split}
\end{equation}
\end{theorem}

It is interesting to remark that the two first terms on the left hand side of \eqref{0.3} is the difference of Fisher's information, while two first terms on
the right hand side is the $2$ times of the difference of entropy.  We mention that in a different framework, some Sobolev estimates for optimal transport
maps have been done in \cite{Figalli1,Figalli2}.

\vskip 2mm
\quad The organization of the paper is as follows. In section 1, we present a construction of the optimal transport map $S$ on the Wiener space $X$,
when the source measure $e^{-W}\mu$ satisfies the Poincar\'e inequality, and target measure $e^{-V}\mu$ is such that the Dirichlet form
${\it E}_V(f,f)=\int_X |\nabla f|_H^2\,e^{-V}d\mu$ is closable; the map $S$ is defined by a $1$-convex function : $S(x)=x+\nabla\psi(x)$ with $\psi\in\D_1^2(X)$.
 In the remainder of the paper, we reverse the source and the target, in order to study the regularity of the inverse map $T$ of $S$. The main task in 
 section 2 is to prove Theorem \ref{thA}: first for a priori estimate, then extended to suitable Sobolev spaces. In section 3, we construct a solution
 to Monge-Amp\`ere equation on the Wiener space: our result (see Theorem \ref{th3.4}) includes two special cases, one studied in \cite{FeyelUstunel2}
 where the source measure is the Wiener measure, another one in \cite{BogachevKolesnikov} where the target measure is the Wiener measure. Besides, 
 we prove that the map $S$ constructed in section 1 admits an inverse map $T$ which is $T(x)=x+\nabla\varphi(x)$
 with $\varphi\in\D_2^2(X)$ (see Theorem 3.5). 

\section{Optimal transport maps on the Wiener space}\label{sect1}
Let $(X,H,\mu)$ be an abstract Wiener space. Consider on $X$ the pseudo-distance $d_H$ defined by
\begin{equation*}
d_H(x,y)=\left\{\begin{array}{ccc}|x-y|_H\ &\hbox{if }& x-y\in H;\\
                                                         +\infty\ &\hbox{otherwise}.& \end{array}\right.
\end{equation*}

Denote by $\P(X)$ the space of probability measures on $X$. For $\nu_1, \nu_2\in\P(X)$, we consider the following Wasserstein distance
\begin{equation*}
W_2^2(\nu_1,\nu_2)=\inf\Bigl\{ \int_{X\times X}d_H(x,y)^2\,\pi(dx,dy);\ \pi\in C(\nu_1,\nu_2)\Bigr\},
\end{equation*}
where $C(\nu_1,\nu_2)$ denotes the totality of probability measures on the product space $X\times X$,
having $\nu_1, \nu_2$ as marginal laws. Note that $W_2(\nu_1,\nu_2)$ could take value $+\infty$. By Talagrand's inequality 
(see for example \cite{Ledoux}),
$W_2^2(\mu, f\mu)\leq 2 \int_X f\log{f}\,d\mu$, that we will denote the latter term by $\Ent_\mu(f)$, we have
\begin{equation}\label{1.1}
W_2(f\mu,g\mu)\leq \sqrt{2}\Bigl(\sqrt{\Ent_\mu(f)}+ \sqrt{\Ent_\mu(g)}\Bigr),
\end{equation}

which is finite, if the measures $f\mu$ and $g\mu$ have finite entropy.  In this situation, it was proven in \cite{FeyelUstunel1}
that there is a unique map $\xi: X\ra H$ such that $x\ra x+\xi(x)$ pushes $f\mu$ to $g\mu$ and 
$W_2(f\mu,g\mu)^2=\int_X |\xi|_H^2\,fd\mu$. However for a general source measure $f\mu$, the construction in \cite{FeyelUstunel1}
is not explicit. For our purpose and the sake of self-contained, we will use the construction in the first part of \cite{FeyelUstunel1}, that is
the usual way when the cost function is strictly convex (see \cite{AGS}, \cite{Villani2}).

\vskip 2mm
\quad Let's introduce some notations in Malliavin calculus (see \cite{Malliavin}, \cite{Fang}). A function $f: X\ra \R$ is called to be cylindrical
if it admits the expression
\begin{equation}\label{1.2}
f(x)=\hat f(e_1(x), \ldots, e_N(x)),\quad \hat f\in C_b^\infty(\R^N), N\geq 1
\end{equation}
where $\{e_1, \ldots, e_N\}$ are elements in dual space $X^*$  of $X$. We denote by $\Cylin(X)$ the space of cylindrical functions on $X$. 
For $f\in\Cylin(X)$
given in \eqref{1.2}, the gradient $\nabla f(x)\in H$ is defined by

\begin{equation}\label{1.3}
\nabla f(x)=\sum_{j=1}^N
\partial_j \hat f(e_1(x), \ldots, e_N(x))\,e_j,
\end{equation}
where $\partial_j$ is ith-partial derivative.
Let $K$ be a separable Hilbert space; a map $F: X\ra K$ is cylindrical if $F$ admits the expression
\begin{equation}\label{1.4}
F=\sum_{i=1}^m f_i k_i,\quad f_i\in\Cylin(X), k_i\in K.
\end{equation}

We denote by $\Cylin(X,K)$ the space of $K$-valued cylindrical functions. For $F\in\Cylin(X,K)$, 
define $\nabla F=\sum_{i=1}^m \nabla f_i\otimes k_i$ which is a $H\otimes K$-valued function. For $h\in H$, we denote
\begin{equation*}
\langle \nabla F,h\rangle = \sum_{i=1}^m\langle \nabla f_i, h\rangle_H\,k_i\in K.
\end{equation*}

In such a way, for any $f\in \Cylin(X)$ and any integer $k\geq 1$, we can define, by induction,
\begin{equation*}
\nabla^k f: X\ra \otimes^k H.
\end{equation*}

Let $p\geq 1$; set
\begin{equation}\label{1.5}
||f||_{D_k^p}^p=\sum_{j=0}^k\int_X ||\nabla^j f(x)||_{\otimes^j H}^p\,d\mu(x),
\end{equation}
here we used the usual convention $\otimes^0 H=\R, \nabla^0f=f$. The Sobolev space $\D_k^p(X)$ is the completion of $\Cylin(X)$ under
the norm defined in \eqref{1.5}. In the same way, the Sobolev space $\D_k^p(X;K)$ of $K$-valued functions is defined. 
\vskip 2mm

Let $V: X\ra\R$ be a measurable function such that $e^{-V}$ is bounded and $\int_X e^{-V}\,d\mu=1$. Consider
\begin{equation}\label{1.6}
{\cal E}_V(F,F)=\int_X ||\nabla F||_{H\otimes K}^2\,e^{-V}d\mu,\quad F\in\Cylin(X,K).
\end{equation}

It is well-known that if 
\begin{equation}\label{1.7}
\int_X |\nabla V|^2\, e^{-V}d\mu<+\infty,
\end{equation}
then the quadratic form \eqref{1.6} is closable over $\Cylin(X,K)$. We will denote by $\D_k^p(X,K;e^{-V}\mu)$ the closure 
of $\Cylin(X,K)$ with respect to the norm defined in \eqref{1.5} replacing $\mu$ by $e^{-V}\mu$.

\vskip 2mm
Let $W\in\D_2^2(X)$ such that $e^{-W}$ is bounded and $\int_X e^{-W}d\mu=1$. Assume that
\begin{equation}\label{1.8}
\nabla^2W\geq -c\,\Id,\quad c\in [0,1[.
\end{equation}

It is known (see \cite{BakryEmery,FeyelUstunel3}) that the condition  \eqref{1.8} implies the following logarithmic Sobolev inequality
\begin{equation}\label{1.9}
(1-c)\int_X {|f|\over ||f||_{L^2(e^{-W}\mu)}}\,e^{-W}d\mu\leq \int_X |\nabla f|^2\,e^{-W}d\mu,\quad f\in\Cylin(X).
\end{equation}

It is also known (see for example \cite{Wang}) that \eqref{1.9} is stronger than Poincar\'e inequality

\begin{equation}\label{1.10}
(1-c)\int_X (f-\E_W(f))^2\, e^{-W}d\mu\leq \int_X |\nabla f|^2\,e^{-W}d\mu,
\end{equation}

where $\E_W$ denotes the integral with respect to the measure $e^{-W}\mu$.

\begin{theorem}\label{th1.1} Under above conditions on $V$ and $W$, there is a  $\psi\in\D_1^2(X,e^{-W}\mu)$
such that $x\ra S(x)=x+\nabla\psi(x)$ is the optimal transport map which pushes $e^{-W}\mu$ to $e^{-V}\mu$; moreover
the inverse map of $S$ is given by $x\ra x+\eta(x)$ with  $\eta\in L^2(X,H;e^{-V}\mu)$.
\end{theorem}

\vskip 2mm
{\bf Proof.} Let $\{e_n;\ n\geq 1\}\subset X^*$ be an orthonormal basis of $H$ and set 
\begin{equation*}
H_n=\hbox{spann}\{e_1, \ldots, e_n\}
\end{equation*}
the vector space spanned by $e_1, \ldots, e_n$, endowed with the induced norm of $H$. Let $\gamma_n$ be the standard Gaussian measure
on $H_n$. Denote
\begin{equation*}
\pi_n(x)=\sum_{j=1}^n e_j(x)\,e_j.
\end{equation*}
Then $\pi_n$ sends the Wiener measure $\mu$ to $\gamma_n$. Let $\F_n$ be the sub $\sigma$-field on $X$ generated by $\pi_n$, 
and  $\E(\ |\F_n)$ be the conditional expectation with respect to $\mu$ and to $\F_n$. Then we can write down

\begin{equation}\label{1.11}
\E(e^{-W}|\F_n)=e^{-W_n}\circ\pi_n,\quad \E(e^{-V}|\F_n)=e^{-V_n}\circ\pi_n.
\end{equation}

Note that for any $f\in L^1(H_n, \gamma_n)$, 
\begin{equation*}
\int_X f\circ\pi_n e^{-W}\,d\mu=\int_X f\circ\pi_n\,\E(e^{-W}|\F_n)\,d\mu=\int_{H_n} fe^{-W_n}d\gamma_n.
\end{equation*}

Applying \eqref{1.10} to $f\circ\pi_n$ yields 
\begin{equation}\label{1.12}
(1-c)\int_{H_n} \Bigl(f-\int_{H_n}fe^{-W_n}d\gamma_n\Bigr)^2\,e^{-W_n}d\gamma_n
\leq \int_{H_n}|\nabla f|^2 e^{-W_n}d\gamma_n,\quad f\in C_b^1(H_n).
\end{equation}

By Kantorovich dual representation theorem (see \cite{Villani2}), we have $W_2^2(e^{-W_n}\gamma_n,e^{-V_n}\gamma_n)=\sup_{(\psi,\varphi)\in\Phi_c}J(\psi,\varphi)$,
where
\begin{equation*}
\Phi_c=\bigl\{(\psi,\varphi)\in L^1(e^{-W_n}\gamma_n)\times  L^1(e^{-V_n}\gamma_n); \psi(x)+\varphi(y)\leq |x-y|_{H_n}^2\bigr\},
\end{equation*}
and
\begin{equation*}
J(\psi,\varphi)=\int_{H_n}\psi(x) e^{-W_n}d\gamma_n+\int_{H_n}\varphi(y)\,e^{-V_n}d\gamma_n.
\end{equation*}
We know there exists a couple of functions $(\psi_n,\varphi_n)$ in $\Phi_c$, which can be chosen to be concave,
such that $W_2^2(e^{-W_n}\gamma_n,e^{-V_n}\gamma_n)=J(\psi_n,\varphi_n)$.
Let $\Gamma_0^n\in C(e^{-W_n}\gamma_n,e^{-V_n}\gamma_n)$ be an optimal coupling, that is,
\begin{equation*}
\int_{H_n\times H_n}|x-y|_{H_n}^2\,d\Gamma_0^n(x,y)=W_2^2(e^{-W_n}\gamma_n,e^{-V_n}\gamma_n).
\end{equation*}

Then it holds true,

\begin{equation}\label{1.13}
|x-y|_{H_n}^2\geq \psi_n(x)+\varphi_n(y),\quad (x,y)\in H_n\times H_n,
\end{equation}
and under $\Gamma_0^n$:

\begin{equation}\label{1.14}
|x-y|_{H_n}^2= \psi_n(x)+\varphi_n(y).
\end{equation}

Combining \eqref{1.13} and \eqref{1.14}, $\Gamma_0^n$ is supported by the graph of $x\ra x-{1\over 2}\nabla\psi_n(x)$ so that
\begin{equation*}
{1\over 4}\int_{H_n}|\nabla\psi_n|^2\,e^{-W_n}d\gamma_n=W_2^2(e^{-W_n}\gamma_n,e^{-V_n}\gamma_n).
\end{equation*}

As in \cite{FeyelUstunel1}, the sequence $\{W_2^2(e^{-W_n}\gamma_n,e^{-V_n}\gamma_n); n\geq 1\}$ is increasing, and converges to
$W_2^2(e^{-W}\mu,e^{-V}\mu)$. Now by \eqref{1.12}, changing $\psi_n$ to $\psi_n-\int_{H_n}\psi_n e^{-W_n}d\gamma_n$,
then $\psi_n\in\D_1^2(e^{-W_n}\gamma_n)$ and
\begin{equation*}
||\psi_n||_{\D_1^2(e^{-W_n}\gamma_n)}^2\leq 2\int_{H_n}|\nabla\psi_n|^2\,e^{-W_n}d\gamma_n.
\end{equation*}
According to \eqref{1.1}, we get that $\sup_{n\geq 1}||\psi_n||_{\D_1^2(e^{-W_n}\gamma_n)}^2<+\infty$.  Now consider $\tilde\psi_n=\psi_n\circ\pi_n$,
$\tilde\varphi_n=\varphi_n\circ\pi_n$. Then
\begin{equation}\label{1.15}
\sup_{n\geq 1}||\tilde\psi_n||_{\D_1^2(e^{-W}\mu)}<+\infty.
\end{equation}

As in \cite{FeyelUstunel1}, define $F_n(x,y)=d_H(x,y)^2-\tilde\psi_n(x)-\tilde\varphi_n(y)$, which is non negative according to \eqref{1.13}. Let
$\Gamma_0$ be an optimal coupling between $e^{-W}\mu$ and $e^{-V}\mu$. We have
\begin{equation}\label{1.16}
\begin{split}
\int_{X\times X}F_n(x,y)\Gamma_0(dx,dy)&=W_2^2(e^{-W}\mu,e^{-V}\mu)-\int_X \tilde\psi_n(x)e^{-W}d\mu-\int_X\tilde\varphi_n(y)\,e^{-V}d\mu\\
&=W_2^2(e^{-W}\mu,e^{-V}\mu)-\int_{H_n} \psi_n(x)e^{-W_n}d\gamma_n-\int_{H_n}\varphi_n(y)\,e^{-V_n}d\gamma_n\\
&=W_2^2(e^{-W}\mu,e^{-V}\mu)-W_2^2(e^{-W_n}\gamma_n,e^{-V_n}\gamma_n)
\end{split}
\end{equation}
which tends to $0$ as $n\ra +\infty$. 
Now returning to \eqref{1.15}, by Banach-Saks theorem, up to a subsequence, the Cesaro mean 
${1\over n}\sum_{j=1}^n \tilde\psi_j$ converges to $\hat\psi$ in $D_1^2(e^{-W}\mu)$. Therefore
\begin{equation*}
{1\over n}\sum_{j=1}^n \tilde\varphi_n(y)=d_H^2(x,y)-{1\over n}\sum_{j=1}^n \tilde\psi_j(x)-{1\over n}\sum_{j=1}^nF_j(x,y)
\end{equation*}
which converges in $L^1$ to $\hat\varphi(y)=d_H^2(x,y)-\hat\psi(x)$.  Now define
\begin{equation*}
\psi=\varliminf_{n\ra +\infty} {1\over n}\sum_{j=1}^n \tilde\psi_j,\quad \varphi=\varliminf_{n\ra +\infty} {1\over n}\sum_{j=1}^n \tilde\varphi_j.
\end{equation*}

Then $\psi=\hat\psi$ for $e^{-W}\mu$ almost all,  $\varphi=\hat\varphi$ for $e^{-V}\mu$ almost all, and by \eqref{1.13}, it holds that
\begin{equation}\label{1.17}
\psi(x)+\varphi(y)\leq d_H^2(x,y),\quad(x,y)\in X\times X.
\end{equation}

Also by above construction, under $\Gamma_0$
\begin{equation}\label{1.18}
\psi(x)+\varphi(y)= d_H^2(x,y).
\end{equation}

Denote by $\Theta_0$  the subset of $(x,y)$ satisfying \eqref{1.18}. On the other hand, the fact that $\psi\in\D_1^2(e^{-W}\mu)$ implies
that for any $h\in H$, there is a full measure subset 
$\Omega_h\subset X$ such that for $x\in \Omega_h$, there is a sequence $\ee_j\downarrow 0$ such that
\begin{equation*}
\langle \nabla\psi(x),h\rangle_H=\lim_{j\ra+\infty} {\psi(x+\ee_j h)-\psi(x)\over\ee_j}.
\end{equation*}

Let $D$ be a countable dense subset of $H$. Then there exists a full measure subset $\Omega$ such that for 
each $x\in\Omega$, for any $h\in D$, there is a sequence $\ee_j\downarrow 0$ such that 
\begin{equation*}
\langle \nabla\psi(x),h\rangle_H=\lim_{j\ra+\infty} {\psi(x+\ee_j h)-\psi(x)\over\ee_j}.
\end{equation*}

Set $\Theta=(\Omega\times X)\cap \Theta_0$. Then $\Gamma_0(\Theta)=1$. For each couple $(x,y)\in\Theta$, we have
$\psi(x)+\varphi(y)=d_H^2(x,y)$ and $\psi(x+\ee_jh)+\varphi(y)\leq d_H^2(x+\ee_j h,y)$. Because $x-y\in H ~~ \Gamma_0-$a.a. it follows that
\begin{equation*}
\psi(x+\ee_jh)-\psi(x)\leq 2\ee_j \langle h, x-y\rangle_H + \ee_j^2|h|_H^2.
\end{equation*}
Therefore $\langle \nabla\psi(x),h\rangle_H \leq 2\langle x-y,h\rangle_H$ for any $h\in D$. From which we deduce that
\begin{equation}\label{1.19}
y=x-{1\over 2}\nabla\psi(x),
\end{equation}
and $\Gamma_0$ is supported by the graph of $x\ra S(x)=x-{1\over 2}\nabla\psi(x)$. Replacing $-{1\over 2}\psi$ by $\psi$, we get 
the statement of the first part of the theorem. For the second part, we refer to section 4 in \cite{FeyelUstunel1}. \fin

\vskip 2mm
For later use, we will emphaze  that   the above constructed whole sequence
\begin{equation}\label{1.20}
 \tilde\varphi_n\ra \varphi\ \hbox{in } L^1(e^{-V}\mu).
 \end{equation}

In fact, if $\tilde\psi$ is another cluster point of $\{\tilde\psi_n; n\geq 1\}$ for the weak topology of $\D_1^2(e^{-W}\mu)$, then 
under the optimal plan $\Gamma_0$, the relation \eqref{1.19} holds for $\tilde\psi$.  Therefore $\nabla\psi=\nabla\tilde\psi$ almost 
everywhere for $e^{-W}\mu$; it follows that $\psi=\tilde\psi$, since $\int_X \psi e^{-W}d\mu=\int_X\tilde\psi\,e^{-W}d\mu=0$.
Now note that 
\begin{equation*}
\begin{split}
\int_X |\nabla\tilde\psi_n|_H^2 e^{-W}d\mu&=\int_{H_n}|\nabla\psi_n|_{H_n}^2 e^{-W_n}d\gamma_n
=W_2^2(e^{-W_n}\gamma_n, e^{-V_n}\gamma_n)\\
&\ra W_2^2(e^{-W}\mu,e^{-V}\mu)=\int_X |\nabla\psi|_H^2 e^{-W}d\mu.
\end{split}
\end{equation*}
Combining these two points, we see that $\tilde\psi_n$ converges to $\psi$ in $\D_1^2(e^{-W}\mu)$. By \eqref{1.16}, the sequence
$\tilde\varphi_n$ converges to $\varphi$ in $L^1(e^{-V}\mu)$. \fin

\section{Variation of optimal transport maps in Sobolev spaces}
\subsection{A priori estimates}

\quad Consider a probability measure $d\mu=e^{-\alpha(x)}\,dx$ on the Euclidean space $(\R^d, |\cdot|)$,
where $\alpha :\mathbb{R}^d\rightarrow \mathbb{R}$ is smooth. Let $h, f$ be two positive 
functions on $\R^d$ such that $\int_{\R^d} h\,d\mu=\int_{\R^d} f\,d\mu=1$. Under some smooth conditions on $h$ and $f$ 
(see \cite{Caffarelli, Kolesnikov1} or p. 561 in \cite{Villani3}), there exists a smooth convex function $\Phi:\R^d\ra\R$ such that $\nabla\Phi: \R^d\ra\R^d$  is a
diffeomorphism which pushes
$h\mu$ forwards to $f\mu$:  $(\nabla\Phi)_\#(h\mu)=f\mu$ and

\begin{equation}\label{2.1}
W_2^2(h\mu,f\mu)=\int_{\R^d}|x-\nabla\Phi(x)|^2\,h(x)d\mu(x),
\end{equation}
where $W_2(h\mu,f\mu)$ denotes the Wasserstein distance between the probability measures $h\mu$ and $f\mu$, which is defined by
\begin{equation*}
W_2^2(h\mu,f\mu)=\inf\Bigl\{\int_{\R^d\times\R^d}|x-y|^2\,d\pi(x,y);\quad \pi\in C(h\mu,f\mu)\Bigr\},
\end{equation*}

the set $C(h\mu,f\mu)$ being the totality of probability measures on the product space $\R^d\times\R^d$ such that $h\mu$ and $f\mu$ are marginals. 

By formula of change of variables, $\nabla\Phi$ satisfies the following
Monge-Amp\`ere equation

\begin{equation}\label{2.2}
f(\nabla\Phi)e^{-\alpha(\nabla\Phi)}\,\det (\nabla^2\Phi)=he^{-\alpha}.
\end{equation}

Now consider two couples of positive functions $(h_1,f_1)$ and $(h_2,f_2)$ satisfying same conditions as $(h,f)$. Let $\Phi_1$ and $\Phi_2$
 be the associated functions. Then we have

\begin{equation}\label{2.3}
f_1(\nabla\Phi_1)e^{-\alpha(\nabla\Phi_1)}\det (\nabla^2\Phi_1)=h_1e^{-\alpha},
\end{equation}

\begin{equation}\label{2.4}
f_2(\nabla\Phi_2)e^{-\alpha(\nabla\Phi_2)}\det (\nabla^2\Phi_2)=h_2e^{-\alpha}.
\end{equation}

Let $S_2$ be the inverse map of $\nabla\Phi_2$, that is,  $\dis \nabla\Phi_2(S_2(x))=x$ on $\R^d$; then we have
\begin{equation*}
\nabla^2\Phi_2(S_2(x))\,\nabla S_2(x)=\Id,\ \hbox{or}\quad \nabla S_2(x)=(\nabla^2\Phi_2)^{-1}(S_2(x)).
\end{equation*}

Acting on the right by $S_2$ the two hand sides of \eqref{2.3}, as well as of \eqref{2.4}, we get

\begin{equation}\label{2.5}
f_1(\nabla\Phi_1(S_2))e^{-\alpha(\nabla\Phi_1(S_2))}\det (\nabla^2\Phi_1(S_2))=h_1(S_2)e^{-\alpha(S_2)},
\end{equation}

\begin{equation}\label{2.6}
f_2\,e^{-\alpha}\,\det (\nabla^2\Phi_2(S_2))=h_2(S_2)e^{-\alpha(S_2)}.
\end{equation}

It follows that

\begin{equation*}
{f_1\over f_2}\cdot {f_1(\nabla\Phi_1(S_2))e^{-\alpha(\nabla\Phi_1(S_2))}\over f_1 e^{-\alpha}}\cdot \det\Bigl[(\nabla^2\Phi_1)(\nabla^2\Phi_2)^{-1}\Bigr](S_2)
={h_1(S_2)\over h_2(S_2)}.
\end{equation*}

Taking the logarithm on the two sides yields
\begin{equation}\label{2.7}
\begin{split}
\log({f_1\over f_2})+& \log(f_1e^{-\alpha})(\nabla\Phi_1(S_2))-\log(
f_1e^{-\alpha})\\
&\hskip 6mm +\log \det\Bigl[(\nabla^2\Phi_1)(\nabla^2\Phi_2)^{-1}\Bigr](S_2)
=\log({h_1\over h_2})(S_2).
\end{split}
\end{equation}

Integrating the two sides of \eqref{2.7} with respect to the measure $f_2\mu$, we get
\begin{equation}\label{2.8}
\begin{split}
\int_{\R^d}\log({h_1\over h_2})(S_2)\, f_2 d\mu -\int_{\R^d}\log({f_1\over f_2})\,f_2 d\mu
=&\int_{\R^d} \log \det\Bigl[(\nabla^2\Phi_1)(\nabla^2\Phi_2)^{-1}\Bigr](S_2)\,f_2 d\mu\\
&\hskip-6mm +\int_{\R^d} \Bigl[ \log(f_1e^{-\alpha})(\nabla\Phi_1(S_2))-\log(
f_1e^{-\alpha})
\Bigr]\, f_2 d\mu.
\end{split}
\end{equation}

By Taylor formula up to order 2,
\begin{equation}\label{2.9}
\begin{split}
& \log(f_1e^{-\alpha})(\nabla\Phi_1(S_2))-\log(
f_1e^{-\alpha})
=\langle \nabla \log(f_1e^{-\alpha}), \nabla\Phi_1(S_2(x))-x\rangle\\
&+\int_0^1(1-t)\Bigl[\nabla^2\log(f_1e^{-\alpha})((1-t)x+t\nabla\Phi_1(S_2(x))\Bigr]\cdot (\nabla\Phi_1(S_2(x))-x)^2\,dt.
\end{split}
\end{equation}

We have
\begin{equation*}
\begin{split}
&\int_{\R^d}\langle \nabla\log(f_1e^{-\alpha}),\nabla\Phi_1(S_2(x))-x\rangle\,f_2d\mu\\
&=\int_{\R^d}\langle\nabla(f_1e^{-\alpha}),\nabla\Phi_1(S_2(x))-x\rangle \,{f_2\over f_1}\,dx.
\end{split}
\end{equation*}

By integration by parts, this last term goes to

\begin{equation*}
\begin{split}
&-\int_{\R^d} f_1e^{-\alpha}\,\div\Bigl(\nabla\Phi_1(S_2(x))-x\Bigr)\,{f_2\over f_1}\,dx
-\int_{\R^d} f_1e^{-\alpha}\langle\nabla\Phi_1(S_2(x))-x, \nabla({f_2\over f_1})\rangle\,dx\\
&=-\int_{\R^d}  \div\Bigl(\nabla\Phi_1(S_2(x))-x\Bigr)\,f_2 d\mu
-\int_{\R^d}  \langle\nabla\Phi_1(S_2(x))-x, \nabla(\log{f_2\over f_1})\rangle\,f_2d\mu.
\end{split}
\end{equation*}

Note that $\dis\nabla\Bigl[(\nabla\Phi_1)(S_2)\Bigr]
=\nabla^2\Phi_1(S_2)\,\nabla S_2=\nabla^2\Phi_1(S_2)\cdot (\nabla^2\Phi_2)^{-1}(S_2)$, and 

\begin{equation*}
\div\Bigl(\nabla\Phi_1(S_2(x))-x\Bigr)=\hbox{Trace}\Bigl[\nabla^2\Phi_1(S_2)\cdot (\nabla^2\Phi_2)^{-1}(S_2)-\Id\Bigr].
\end{equation*}

Combining above computations yields
\begin{equation}\label{2.10}
\begin{split}
&\int_{\R^d}\langle \nabla\log(f_1e^{-\alpha}), \nabla\Phi_1(S_2(x))-x\rangle\,f_2d\mu\\
=&-\int_{\R^d} \hbox{Trace}\Bigl[\nabla^2\Phi_1(S_2)\cdot (\nabla^2\Phi_2)^{-1}(S_2)-\Id\Bigr]\,f_2d\mu\\
&-\int_{\R^d}  \langle\nabla\Phi_1(S_2(x))-x, \nabla(\log{f_2\over f_1})\rangle\,f_2d\mu.
\end{split}
\end{equation}

For a matrix $A$ on $\R^d$, the Fredholm-Carleman determinant $\det_2(A)$ is defined by
\begin{equation*}
\det_2(A)=e^{\hbox{Trace}(\Id-A)}\,\det(A).
\end{equation*}

It is easy to check that if $A$ is symmetric positive, then $\dis 0\leq\det_2(A)\leq 1$. We have
\begin{equation*}
 \hbox{Trace}\Bigl((\nabla^2\Phi_1)(\nabla^2\Phi_2)^{-1}\Bigr)
 = \hbox{Trace}\Bigl((\nabla^2\Phi_2)^{-1/2}\,\nabla^2\Phi_1\,(\nabla^2\Phi_2)^{-1/2}\Bigr),
 \end{equation*}
and 
\begin{equation*}
 \det\Bigl((\nabla^2\Phi_1)(\nabla^2\Phi_2)^{-1}\Bigr)
 = \det\Bigl((\nabla^2\Phi_2)^{-1/2}\,\nabla^2\Phi_1\,(\nabla^2\Phi_2)^{-1/2}\Bigr).
 \end{equation*}
 
 Therefore
 \begin{equation}\label{2.11}
 \log\det_2\Bigl((\nabla^2\Phi_1)(\nabla^2\Phi_2)^{-1}\Bigr)
 = \log\det_2\Bigl((\nabla^2\Phi_2)^{-1/2}\,\nabla^2\Phi_1\,(\nabla^2\Phi_2)^{-1/2}\Bigr)\leq 0.
 \end{equation}

Now combining \eqref{2.8}, \eqref{2.9} and \eqref{2.10}, we get the following result.

\begin{theorem}\label{th2.1}  Let $\alpha\in C^\infty(\R^d)$ and $d\mu=e^{-\alpha}dx$ be a probability measure on $\R^d$. Then
\begin{equation}\label{2.12}
\begin{split}
&\Ent_{h_1\mu}\bigl({h_2\over h_1}\bigr)-\Ent_{f_1\mu}\bigl({f_2\over f_1}\bigr)
=\int_{\R^d}\langle\nabla\Phi_1-\nabla\Phi_2, \nabla(\log{f_2\over f_1})(\nabla\Phi_2)\rangle\, h_2d\mu\\
&-\int_{\R^d} \log\det_2\Bigl((\nabla^2\Phi_2)^{-1/2}\,\nabla^2\Phi_1\,(\nabla^2\Phi_2)^{-1/2}\Bigr)\,h_2d\mu\\
&+\int_0^1(1-t)dt\int_{\R^d}\Bigl[-\nabla^2\log(f_1e^{-\alpha})((1-t)\nabla\Phi_2+t\nabla\Phi_1)\Bigr]\cdot (\nabla\Phi_1-\nabla\Phi_2)^2\,h_2d\mu.
\end{split}
\end{equation}
\end{theorem}

\begin{corollary}\label{th2.2}
Suppose that
\begin{equation}\label{2.13}
\nabla^2\bigl(-\log(f_1e^{-\alpha})\bigr)\geq c\,\Id,\quad c>0.
\end{equation}
Then
\begin{equation}\label{2.14}
\begin{split}
\int_{\R^d}|\nabla\Phi_1-\nabla\Phi_2|^2\,h_2d\mu
&\leq {4\over c}\Bigl(\Ent_{h_1\mu}\bigl({h_2\over h_1}\bigr)-\Ent_{f_1\mu}\bigl({f_2\over f_1}\bigr)\Bigr)\\
&+{4\over c^2} \int_{\R^d}|\nabla\log{f_2\over f_1}|^2\, f_2d\mu.
\end{split}
\end{equation}
If moreover $f_1=f_2$, then it holds more precisely
\begin{equation*}
{c\over 2}\,\int_{\R^d}|\nabla\Phi_1-\nabla\Phi_2|^2\,h_2d\mu
\leq \Ent_{h_1\mu}\bigl({h_2\over h_1}\bigr).
\end{equation*}
\end{corollary}

\vskip 2mm
{\bf Proof.} Note that
\begin{equation*}
\begin{split}
\Bigl|\int_{\R^d}\langle\nabla\Phi_1-\nabla\Phi_2, \nabla(\log{f_2\over f_1})(\nabla\Phi_2)\rangle\, h_2d\mu\Bigr|
&\leq \Bigl(\int_{\R^d}|\nabla\Phi_1-\nabla\Phi_2|^2\, h_2d\mu\Bigr)^{1/2}\Bigl(\int_{\R^d}|\nabla\log{f_2\over f_1}|^2 \,f_2d\mu\Bigr)^{1/2}\\
&\leq {c\over 4}\int_{\R^d}|\nabla\Phi_1-\nabla\Phi_2|^2\, h_2d\mu+ {1\over c}\int_{\R^d}|\nabla\log{f_2\over f_1}|^2\,f_2d\mu.
\end{split}
\end{equation*}

Under  condition \eqref{2.13}, the last term in \eqref{2.12} is bounded from below by
\begin{equation*}
{c\over 2}\int_{\R^d}|\nabla\Phi_1-\nabla\Phi_2|^2\,h_2d\mu.
\end{equation*}

Now according to \eqref{2.12}, we get the result from \eqref{2.14}. \fin

\vskip 2mm
\quad In what follows, we will consider the standard Gaussian measure $\gamma$ as the reference measure on $\R^d$. Let $e^{-V}$ and $e^{-W}$ be two 
density functions with respect to $\gamma$, that is, $\int_{\R^d} e^{-V}d\gamma=\int_{\R^d} e^{-W}d\gamma=1$. Let $\Phi$ be a smooth convex function such that
$\nabla\Phi$ pushes $e^{-V}\gamma$ forward to $e^{-W}\gamma$, that is,  

\begin{equation*}
\int_{\R^d}F(\nabla\Phi)\,e^{-V}d\gamma=\int_{\R^d} F\,e^{-W}d\gamma.
\end{equation*}

Let $a\in\R^d$; then
\begin{equation*}
\int_{\R^d}F(\nabla\Phi(x+a))e^{-V(x+a)}e^{-\langle x,a\rangle-{1\over2}|a|^2}\,d\gamma
=\int_{\R^d} F(\nabla\Phi) e^{-V}\,d\gamma.
\end{equation*}

Denote by $\tau_a$ the translation by $a$, and $M_a(x)=e^{-\langle x,a\rangle-{1\over2}|a|^2}$, then the above relations imply that
\begin{equation*}
\nabla(\tau_a\Phi)_\#: e^{-\tau_aV}M_a\gamma\ra e^{-W}\gamma.
\end{equation*}

Let $h_1=e^{-\tau_aV}M_a, h_2=e^{-V}$ . Then $\Ent_{h_1\mu}\bigl({h_2\over h_1}\bigr)=\int_{\R^d} (\tau_aV-V+\langle x,a\rangle+{1\over2}|a|^2) e^{-V}d\gamma$.
Applying Theorem \ref{th2.1} , we get
\begin{equation*}
\begin{split}
&\int_{\R^d} (\tau_aV-V+\langle x,a\rangle+{1\over2}|a|^2) e^{-V}d\gamma\\
=&-\int_{\R^d} \log\det_2\Bigl[(\nabla^2\Phi)^{-1/2}\,\nabla^2(\tau_a\Phi)\,(\nabla^2\Phi)^{-1/2}\Bigr]\,e^{-V}d\gamma\\
&+ \int_0^1 (1-t)dt\int_{\R^d}\Bigl[(\Id+\nabla^2W)(\Lambda(t,x,a))\Bigr]\cdot (\nabla\Phi(x)-\nabla\Phi(x+a))^2 e^{-V}d\gamma,
\end{split}
\end{equation*}
where $\Lambda(t,x,a)=(1-t)\nabla\Phi(x)+t\nabla\Phi(x+a)$. Note that as $a\ra0$, $\Lambda(t,x,a)\ra \nabla\Phi(x)$.

\quad Replacing $a$ by $-a$, and summing respectively the two hand sides of these equalities, we get
\begin{equation}\label{A1}
\begin{split}
&\int_{\R^d} \bigl(V(x+a)+V(x-a)-2V(x)+|a|^2\bigr)\, e^{-V}d\gamma = J(a)+J(-a)\\
&+ \int_0^1 (1-t)dt\int_{\R^d}\Bigl[(\Id+\nabla^2W)(\Lambda(t,x,a))\Bigr]\cdot (\nabla\Phi(x)-\nabla\Phi(x+a))^2 e^{-V}d\gamma\\
 &+\int_0^1 (1-t)dt\int_{\R^d}\Bigl[(\Id+\nabla^2W)(\Lambda(t,x,-a))\Bigr]\cdot (\nabla\Phi(x)-\nabla\Phi(x-a))^2 e^{-V}d\gamma,
\end{split}
\end{equation}
where 
\begin{equation*}
J(a)=-\int_{\R^d} \log\det_2\Bigl[(\nabla^2\Phi)^{-1/2}\,\nabla^2(\tau_a\Phi)\,(\nabla^2\Phi)^{-1/2}\Bigr]\,e^{-V}d\gamma.
\end{equation*}
By explicit formula in Lemma \ref{th4.1} in appendice, and write $\nabla\Phi(x)=x+\nabla\varphi(x)$, we have
\begin{equation*}
\begin{split}
&{1\over\ee^2}J(\ee a)=\int_0^1(1-t)dt\int_{\R^d}||(I+(1-t)\nabla^2\varphi+t\nabla^2\varphi(x+\ee a))^{-1/2}\\
&\hskip 15mm \ee^{-1}\Bigl(\nabla^2\varphi(x+\ee a)-\nabla^2\varphi(x)\Bigr)
(I+(1-t)\nabla^2\varphi+t\nabla^2\varphi(x+\ee a))^{-1/2}||_{HS}^2e^{-V}d\gamma.
\end{split}
\end{equation*}

So that, by Fatou lemma
\begin{equation}\label{A2}
\varliminf_{\ee\ra0}{J(\ee a)\over\ee^2}\geq{1\over2}\int_{\R^d}||(I+\nabla^2\varphi)^{-1/2}\,D_a\nabla^2\varphi(x)\,(I+\nabla^2\varphi)^{-1/2}||_{HS}^2\,e^{-V}d\gamma.
\end{equation}

Now replacing $a$ by $\ee a$ and dividing by $\ee^2$ the two hand sides of \eqref{A1}, letting $\ee\ra0$ yields 
\begin{equation}\label{2.15}
\begin{split}
\int_{\R^d} \Bigl[ D_a^2V +|a|^2\Bigr]\, e^{-V}d\gamma&\geq\int_{\R^d}||(I+\nabla^2\varphi)^{-1/2}\,D_a\nabla^2\varphi(x)\,(I+\nabla^2\varphi)^{-1/2}||_{HS}^2\,e^{-V}d\gamma\\
&+ \int_{\R^d}  (\Id+\nabla^2 W)(\nabla\Phi)\,(D_a\nabla\Phi, D_a\nabla\Phi)\, e^{-V}d\gamma\\
&=\int_{\R^d}||(I+\nabla^2\varphi)^{-1/2}\,D_a\nabla^2\varphi(x)\,(I+\nabla^2\varphi)^{-1/2}||_{HS}^2\,e^{-V}d\gamma\\
&+\int_{\R^d} |D_a\nabla\Phi|^2 e^{-V}d\gamma +\int_{\R^d}(\nabla^2W)(\nabla\Phi)(D_a\nabla\Phi,D_a\nabla\Phi)\,e^{-V}d\gamma.
\end{split}
\end{equation}

By integration by parts,
\begin{equation*}
\int_{\R^d} D_a^2 V\,e^{-V}d\gamma
=\int_{\R^d} (D_aV)^2 e^{-V}d\gamma + \int_{\R^d} D_aV\,\langle a,x\rangle\,e^{-V}d\gamma.
\end{equation*}

Using \eqref{2.15} and $|D_a\nabla\Phi|^2=|a|^2 + 2\langle a, D_a\nabla\varphi\rangle +|D_a\nabla\varphi|^2$, we get
\begin{equation*}
\begin{split}
&\int_{\R^d} (D_aV)^2 e^{-V}d\gamma + \int_{\R^d} D_aV\,\langle a,x\rangle\,e^{-V}d\gamma\\
& \geq\int_{\R^d}||(I+\nabla^2\varphi)^{-1/2}\,D_a\nabla^2\varphi(x)\,(I+\nabla^2\varphi)^{-1/2}||_{HS}^2\,e^{-V}d\gamma\\
 &+2\int_{\R^d} \langle a,D_a\nabla\varphi\rangle\,e^{-V}d\gamma
  +\int_{\R^d}|D_a\nabla\varphi|^2\,e^{-V}d\gamma
 +\int_{\R^d}\nabla^2W_{\nabla\Phi}(D_a\nabla\Phi,D_a\nabla\Phi)\,e^{-V}d\gamma.
\end{split}
\end{equation*}

Summing $a$ on an orthonormal basis $ {\cal B}$, it follows

\begin{equation}\label{2.16}
\begin{split}
&\int_{\R^d} |\nabla V|^2 e^{-V}d\gamma + \int_{\R^d} \langle x,\nabla V\rangle\,e^{-V}d\gamma\\
&\hskip -10mm \geq\int_{\R^d}\sum_{a\in{\cal B}}||(I+\nabla^2\varphi)^{-1/2}\,D_a\nabla^2\varphi(x)\,(I+\nabla^2\varphi)^{-1/2}||_{HS}^2\,e^{-V}d\gamma\\
&\hskip -10mm+2\int_{\R^d}\Delta\varphi\,e^{-V}d\gamma 
+\int_{\R^d}||\nabla^2\varphi||_{HS}^2 e^{-V}d\gamma
+\sum_{a\in {\cal B}}\int_{\R^d} \nabla^2W_{\nabla\Phi}(D_a\nabla\Phi,D_a\nabla\Phi)\,e^{-V}d\gamma.
\end{split}
\end{equation}

Let
\begin{equation}\label{2.17}
N_W(\nabla^2\varphi)=\sum_{a\in{\cal B}}\nabla^2W_{\nabla\Phi}(D_a\nabla\varphi,D_a\nabla\varphi).
\end{equation}

Then
\begin{equation*}
\begin{split}
&\sum_{a\in {\cal B}}\int_{\R^d} \nabla^2W_{\nabla\Phi}(D_a\nabla\Phi,D_a\nabla\Phi)\,e^{-V}d\gamma\\
&\hskip -10mm= \int_{\R^d}(\Delta W)(\nabla\Phi)\,e^{-V}d\gamma +2\int_{\R^d}\langle \nabla^2W(\nabla\Phi),\nabla^2\varphi\rangle_{HS}\,e^{-V}d\gamma
+\int_{\R^d} N_W(\nabla^2\varphi)\,e^{-V}d\gamma.
\end{split}
\end{equation*}

This equality, together with \eqref{2.16} yield
\begin{equation}\label{2.18}
\begin{split}
&\int_{\R^d} |\nabla V|^2 e^{-V}d\gamma + \int_{\R^d} \langle x,\nabla V\rangle\,e^{-V}d\gamma\\
&\hskip -10mm\geq\int_{\R^d}\sum_{a\in{\cal B}}||(I+\nabla^2\varphi)^{-1/2}\,D_a\nabla^2\varphi(x)\,(I+\nabla^2\varphi)^{-1/2}||_{HS}^2\,e^{-V}d\gamma\\
&\hskip -10mm + 2\int_{\R^d}\Delta\varphi\,e^{-V}d\gamma +
\int_{\R^d}||\nabla^2\varphi||_{HS}^2 e^{-V}d\gamma+ \int_{\R^d}(\Delta W)(\nabla\Phi)\,e^{-V}d\gamma\\
&\hskip -10mm
+2\int_{\R^d}\langle \nabla^2W(\nabla\Phi),\nabla^2\varphi\rangle_{HS}\,e^{-V}d\gamma+\int_{\R^d} N_W(\nabla^2\varphi)\,e^{-V}d\gamma.
\end{split}
\end{equation}

In order to obtain desired terms, we first use the relation
\begin{equation*}
\int_{\R^d} |x+\nabla\varphi(x)|^2\, e^{-V}d\gamma
=\int_{\R^d} |x|^2\, e^{-W}d\gamma
\end{equation*}

which gives that
\begin{equation*}
2\int_{\R^d}\langle x,\nabla\varphi(x)\rangle\,e^{-V}d\gamma
=\int_{\R^d} |x|^2\, e^{-W}d\gamma-\int_{\R^d} |x|^2\, e^{-V}d\gamma
-\int_{\R^d}|\nabla\varphi(x)|^2\,e^{-V}d\gamma.
\end{equation*}

Let $\L$ be the Ornstein-Uhlenbeck operator: $\L f(x)=\Delta f (x)- \langle x,\nabla f\rangle$. Remark that
\begin{equation*}
\L({1\over2}|x|^2)=d-|x|^2.
\end{equation*}

Then $\int_{\R^d} |x|^2\, e^{-W}d\gamma-\int_{\R^d} |x|^2\, e^{-V}d\gamma
=-\int_{\R^d}\L({1\over2}|x|^2) e^{-W}d\gamma + \int_{\R^d}\L({1\over2}|x|^2) e^{-V}d\gamma$, which is equal to

\begin{equation*}
-\int_{\R^d}\langle x,\nabla W\rangle\,e^{-W}d\gamma +\int_{\R^d}\langle x,\nabla V\rangle\,e^{-V}d\gamma.
\end{equation*}

Therefore
\begin{equation}\label{2.19}
\begin{split}
2 \int_{\R^d}\langle x,\nabla\varphi(x)\rangle\,e^{-V}d\gamma=&-\int_{\R^d}\langle x,\nabla W\rangle\,e^{-W}d\gamma\\
&+\int_{\R^d}\langle x,\nabla V\rangle\,e^{-V}d\gamma -\int_{\R^d}|\nabla\varphi|^2\,e^{-V}d\gamma.
\end{split}
\end{equation}

On the other hand, from Monge-Amp\`ere equation,
\begin{equation*}
e^{-V}=e^{-W(\nabla\Phi)}e^{\L\varphi -{1\over2}|\nabla\varphi|^2}\det_2(\Id+\nabla^2\varphi),
\end{equation*}
we have
\begin{equation*}
-V=-W(\nabla\Phi)+\L\varphi -{1\over2}|\nabla\varphi|^2 +\log\det_2(\Id+\nabla^2\varphi).
\end{equation*}

Integrating the two hand sides with respect to $e^{-V}d\gamma$, we get
\begin{equation}\label{2.20}
\begin{split}
\int_{\R^d}\L\varphi\,e^{-V}d\gamma=&\Ent_\gamma(e^{-V})-\Ent_\gamma(e^{-W})+{1\over2}\int_{\R^d}|\nabla\varphi|^2\,e^{-V}d\gamma\\
&-\int_{\R^d}\log\det_2(\Id+\nabla^2\varphi)\, e^{-V}d\gamma.
\end{split}
\end{equation}

Combining \eqref{2.19} and \eqref{2.20}, we get
\begin{equation*}
\begin{split}
2\int_{\R^d} \Delta\varphi\,e^{-V}d\gamma&=2\int_{\R^d} \L\varphi\,e^{-V}d\gamma+2\int_{\R^d}\langle x,\nabla\varphi\rangle\,e^{-V}d\gamma\\
&=2\Ent_\gamma(e^{-V})-2\Ent_\gamma(e^{-W})-2\int_{\R^d}\log\det_2(\Id+\nabla^2\varphi)\, e^{-V}d\gamma\\
&-\int_{\R^d}\langle x,\nabla W\rangle\,e^{-W}d\gamma +\int_{\R^d}\langle x,\nabla V\rangle\,e^{-V}d\gamma.
\end{split}
\end{equation*}

Replacing $\int_{\R^d} \Delta\varphi\,e^{-V}d\gamma$ in \eqref{2.18} by above expression, we obtain
\begin{equation*}
\begin{split}
&\int_{\R^d}|\nabla V|^2\, e^{-V}d\gamma\geq 2\Ent_\gamma(e^{-V})-2\Ent_\gamma(e^{-W})
-2\int_{\R^d}\log\det_2(\Id+\nabla^2\varphi)\, e^{-V}d\gamma\\
&\hskip5mm+\int_{\R^d}\sum_{a\in{\cal B}}||(I+\nabla^2\varphi)^{-1/2}\,D_a\nabla^2\varphi(x)
\,(I+\nabla^2\varphi)^{-1/2}||_{HS}^2\,e^{-V}d\gamma+\int_{\R^d}||\nabla^2\varphi||_{HS}^2\,e^{-V}d\gamma\\
&\hskip 5mm+\int_{\R^d}\L W\,e^{-W}d\gamma
+2\int_{\R^d}\langle\nabla^2W(\nabla\Phi),\nabla^2\varphi\rangle_{HS}\,e^{-V}d\gamma+\int_{\R^d} N_W(\nabla^2\varphi)\,e^{-V}d\gamma.
\end{split}
\end{equation*}

So we get
\begin{theorem}\label{th2.3} We have 
\begin{equation*}
\begin{split}
&\int_{\R^d}|\nabla V|^2\, e^{-V}d\gamma-\int_{\R^d}|\nabla W|^2\, e^{-W}d\gamma\\
&\hskip -5mm\geq 2\Ent_\gamma(e^{-V})-2\Ent_\gamma(e^{-W})
-2\int_{\R^d}\log\det_2(\Id+\nabla^2\varphi)\, e^{-V}d\gamma\\
&\hskip -5mm+\int_{\R^d}\sum_{a\in{\cal B}}||(I+\nabla^2\varphi)^{-1/2}\,D_a\nabla^2\varphi(x)
\,(I+\nabla^2\varphi)^{-1/2}||_{HS}^2\,e^{-V}d\gamma+\int_{\R^d}||\nabla^2\varphi||_{HS}^2\,e^{-V}d\gamma\\
&\hskip -5mm+2\int_{\R^d}\langle\nabla^2W(\nabla\Phi),\nabla^2\varphi\rangle_{HS}\,e^{-V}d\gamma+\int_{\R^d} N_W(\nabla^2\varphi)\,e^{-V}d\gamma.
\end{split}
\end{equation*}
\end{theorem}

\begin{theorem}\label{th2.4} Assume that $\nabla^2 W\geq -c\,\Id$ with $c\in [0,1[$; then
\begin{equation}\label{2.21}
\begin{split}
&\int_{\R^d}|\nabla V|^2\, e^{-V}d\gamma-\int_{\R^d}|\nabla W|^2\, e^{-W}d\gamma
+{2\over 1-c}\int_{\R^d}||\nabla^2W||_{HS}^2 e^{-W}d\gamma\\
&\geq 2\Ent_\gamma(e^{-V})-2\Ent_\gamma(e^{-W})
+{1-c\over2}\int_{\R^d}||\nabla^2\varphi||_{HS}^2\,e^{-V}d\gamma.
\end{split}
\end{equation}
\end{theorem}

\vskip 2mm
{\bf Proof.} It is sufficient to notice that 
\begin{equation*}
2\int_{\R^d}|\langle\nabla^2W(\nabla\Phi),\nabla^2\varphi\rangle_{HS}|\,e^{-V}d\gamma
\leq  {1-c\over 2}\int_{\R^d}||\nabla^2\varphi||_{HS}^2\,e^{-V}d\gamma +
{2\over 1-c}\int_{\R^d}||\nabla^2W||_{HS}^2\,e^{-W}d\gamma.
\end{equation*}
The inequality \eqref{2.21} follows from Theorem \ref{th2.3}. \fin

\begin{theorem}\label{th2.5} Let $1\leq p<2$. Denote by $||\cdot||_{op}$ the norm of operator, then
\begin{equation}\label{th2.5-1}
||\nabla^3\varphi||_{L^p(e^{-V}\gamma)}^2\leq \Bigl\| ||I+\nabla^2\varphi||_{op}\Bigr\|_{L^{2p\over 2-p}(e^{-V}\gamma)}^2\Bigl(||\nabla V||_{L^2(e^{-V}\gamma)}^2
+{2\over 1-c}||\nabla^2W||_{L^2(e^{-W}\gamma)}^2\Bigr).
\end{equation}
\end{theorem}

\vskip 2mm
{\bf Proof.}
By H\"older inequality
\begin{equation*}
\int_{\R^d}||\nabla^3\varphi||_{HS}^p\,e^{-V}d\gamma
\leq \Bigl(\int_{\R^d}{||\nabla^3\varphi||_{HS}^2\over ||I+\nabla^2\varphi||_{op}^2}\,e^{-V}d\gamma\Bigr)^{p/2}\,
\Bigl(\int_{\R^d}||I+\nabla^2\varphi||_{op}^{2p\over 2-p}\,e^{-V}d\gamma\Bigr)^{2-p\over2}.
\end{equation*}

 By \eqref{2.22} below : 
\begin{equation*}
 {||\nabla^3\varphi||_{HS}^2\over ||I+\nabla^2\varphi||_{op}^2}\leq \sum_{a\in{\cal B}}||(I+\nabla^2\varphi)^{-1/2}\,D_a\nabla^2\varphi(x)\,(I+\nabla^2\varphi)^{-1/2}||_{HS}^2.
\end{equation*}

Remark that $\int_{\R^d}|\nabla W|^2 e^{-W}d\gamma\geq 2\Ent_\gamma(e^{-W})$. Now by Theorem \ref{th2.3}, we get the result.\fin

\vskip 2mm
\quad In what follows, we will compute the variation of optimal transport maps in Sobolev spaces.  Consider
\begin{equation*}
(\nabla\Phi_1)_\#: e^{-V_1}d\gamma\ra e^{-W_1}d\gamma,\quad (\nabla\Phi_2)_\#: e^{-V_2}d\gamma\ra e^{-W_2}d\gamma.
\end{equation*}

We will explore the term $-\log\det_2\Bigl[(\nabla^2\Phi_2)^{-1/2}\nabla^2\Phi_1(\nabla^2\Phi_2)^{-1/2}\Bigr]$ in Theorem \ref{th2.1}.\vskip 2mm
Let $\nabla\Phi_1(x)=x+\nabla\varphi_1(x)$ and $\nabla\Phi_2(x)=x+\nabla\varphi_2(x)$; then
\begin{equation*}
\nabla^2\Phi_1=I+\nabla^2\varphi_1,\quad \nabla^2\Phi_2=I+\nabla^2\varphi_2.
\end{equation*}

\begin{theorem}\label{th2.6}  Let $1\leq p<2$ and
 \begin{equation}\label{2.24}
 M(\nabla^2\varphi_1,\nabla^2\varphi_2)
 =\max\Bigl(\Bigl\|||I+\nabla^2\varphi_1||_{op}\Bigr\|_{L^{2p\over2-p}(e^{-V_2}\gamma)}^2,
  \Bigl\|||I+\nabla^2\varphi_2||_{op}\Bigr\|_{L^{2p\over2-p}(e^{-V_2}\gamma)}^2\Bigr).
 \end{equation}
Assume that $\nabla^2W_1\geq -c\,\Id$ with $c\in [0,1[$. Then we have
\begin{equation}\label{2.25}
\begin{split}
||\nabla^2\varphi_1-\nabla^2\varphi_2||_{L^p(e^{-V_2}\gamma)}^2 
\leq &2M(\nabla^2\varphi_1,\nabla^2\varphi_2) \Bigl[2\int_{\R^d} (V_1-V_2)e^{-V_2}d\gamma \\
&+ {2\over1- c}\int_{\R^d}|\nabla(W_1-W_2)|^2e^{-W_2}d\gamma\Bigr].
\end{split}
\end{equation}
\end{theorem}

\vskip 2mm
{\bf Proof.} Applying Lemma \ref{th4.1} to $B=\nabla^2\varphi_1-\nabla^2\varphi_2$ and $A=I+(1-t)\nabla^2\varphi_2+t\nabla^2\varphi_1$ yields

\begin{equation*}
\begin{split}
&||(I+(1-t)\nabla^2\varphi_2+t\nabla^2\varphi_1)^{-1/2}(\nabla^2\varphi_1-\nabla^2\varphi_2)(I+(1-t)\nabla^2\varphi_2+t\nabla^2\varphi_1)^{-1/2}||_{HS}^2
\\&\hskip 15mm \geq {||\nabla^2\varphi_1-\nabla^2\varphi_2||_{HS}^2\over ||I+(1-t)\nabla^2\varphi_2+t\nabla^2\varphi_1||_{op}^2}.
\end{split}
\end{equation*}

As above, by H\"older inequality, we have
\begin{equation*}
\int_{\R^d} {||\nabla^2\varphi_1-\nabla^2\varphi_2||_{HS}^2\over ||I+(1-t)\nabla^2\varphi_2+t\nabla^2\varphi_1||_{op}^2}\, e^{-V_2}d\gamma
\geq {||\nabla^2\varphi_1-\nabla^2\varphi_2||_{L^p(e^{-V_2}\gamma)}^2 \over
\Bigl\| ||I+(1-t)\nabla^2\varphi_2+t\nabla^2\varphi_1||_{op}\Bigr\|_{L^{2p\over 2-p}(e^{-V_2}\gamma)}^2}.
\end{equation*}

Now by convexity, 
\begin{equation*}
\begin{split}
&\Bigl\| ||I+(1-t)\nabla^2\varphi_2+t\nabla^2\varphi_1||_{op}\Bigr\|_{L^{2p\over 2-p}(e^{-V_2}\gamma)}^2\\
&\hskip -10mm\leq(1-t) \Bigl\| ||I+\nabla^2\varphi_2||_{op}\Bigr\|_{L^{2p\over 2-p}(e^{-V_2}\gamma)}^2+t \Bigl\| ||I+\nabla^2\varphi_1||_{op}\Bigr\|_{L^{2p\over 2-p}(e^{-V_2}\gamma)}^2
\leq M(\nabla^2\varphi_1,\nabla^2\varphi_2).
\end{split}
 \end{equation*}
 
 According to Lemma \ref{th4.2}, we have
\begin{equation}\label{2.26}
\begin{split}
&\int_{\R^d}-\log\det_2\Bigl((\nabla^2\Phi_2)^{-1/2}\,\nabla^2\Phi_1\,(\nabla^2\Phi_2)^{-1/2}\Bigr)\, e^{-V_2}d\gamma\\
&\geq \int_0^1(1-t)dt\int_{\R^d} {||\nabla^2\varphi_1-\nabla^2\varphi_2||_{HS}^2\over ||I+(1-t)\nabla^2\varphi_2+t\nabla^2\varphi_1||_{op}^2}\, e^{-V_2}d\gamma\\
&\geq {1\over 2}{||\nabla^2\varphi_1-\nabla^2\varphi_2||_{L^p(e^{-V_2}\gamma)}^2\over M(\nabla^2\varphi_1,\nabla^2\varphi_2)}.
\end{split}
\end{equation}

 By Cauchy-Schwarz inequality,
\begin{equation*}
\begin{split}
&\Bigl|\int_{\R^d}\langle\nabla\Phi_1-\nabla\Phi_2,\nabla(W_1-W_2)(\nabla\Phi_2)\rangle\,e^{-V_2}d\gamma\Bigr|\\
&\leq \Bigl(\int_{\R^d}|\nabla\Phi_1-\nabla\Phi_2|^2\,e^{-V_2}d\gamma\Bigr)^{1/2}\,\Bigl(\int_{\R^d}|\nabla(W_1-W_2)|^2\,e^{-W_2}d\gamma\Bigr)^{1/2}\\
&\leq \frac{1-c}{4}\int_{\R^d}|\nabla\Phi_1-\nabla\Phi_2|^2\,e^{-V_2}d\gamma+\frac{1}{1-c}\int_{\R^d}|\nabla(W_1-W_2)|^2\,e^{-W_2}d\gamma.
\end{split}
\end{equation*}

Under the hypothesis $\nabla^2W_1\geq -c\Id$ with $c<1$, the inequality \eqref{2.14} implies 
\begin{eqnarray*}
& ~&\int_{\R^d}|\nabla\Phi_1-\nabla\Phi_2|^2\,e^{-V_2}d\gamma
\leq {4\over 1-c}\int_{\R^d}(V_1-V_2)e^{-V_2}d\gamma + {4\over (1-c)^2}\int_{\R^d}|\nabla(W_1-W_2)|^2e^{-W_2}d\gamma,
\end{eqnarray*}

so that 
\begin{equation*}
\begin{split}
&\Bigl|\int_{\R^d}\langle\nabla\Phi_1-\nabla\Phi_2,\nabla(W_1-W_2)(\nabla\Phi_2)\rangle\,e^{-V_2}d\gamma\Bigr|\\
&\leq \int_{\R^d}(V_1-V_2)e^{-V_2}d\gamma+{2\over 1-c}\int_{\R^d}|\nabla(W_1-W_2)|^2e^{-W_2}d\gamma.
\end{split}
\end{equation*}

Now combinig  \eqref{2.12} and \eqref{2.26}, we conclude \eqref{2.25}. \fin

\subsection{Extension to Sobolev spaces}

\quad In this subsection, we will assume that $V\in\D_1^2(\R^d,\gamma), W\in\D_2^2(\R^d,\gamma)$ and there exist constants  $\delta_2>0$ and $c\in [0,1[$ such that
\begin{equation}\label{2.28}
e^{-V}\leq\delta_2,\quad  e^{-W}\leq\delta_2\quad\hbox{and}\ \nabla^2W\geq -c\,\Id.
\end{equation}

It turns out that  $V$ and $W$ are bounded from below. Consider the Ornstein-Uhlenbeck semi-group $P_\ee$
\begin{equation*}
P_\ee f(x)=\int_{\R^d} f(e^{-\ee}x+\sqrt{1-e^{2\ee}}\,y)\,d\gamma(y).
\end{equation*}

If $f\in\D_2^2(\R^d,\gamma)$, then
\begin{equation*}
\nabla P_\ee f(x)=e^{-\ee}\int_{\R^d} \nabla f(e^{-\ee}x+\sqrt{1-e^{2\ee}}\,y)\,d\gamma(y),
\end{equation*}
and
\begin{equation*}
\nabla^2 P_\ee f(x)=e^{-2\ee}\int_{\R^d} \nabla^2 f(e^{-\ee}x+\sqrt{1-e^{2\ee}}\,y)\,d\gamma(y).
\end{equation*}

It follows that $||\nabla P_\ee f||_{L^2(\gamma)}\leq ||\nabla f||_{L^2(\gamma)}$ and $||\nabla^2 P_\ee f||_{L^2(\gamma)}\leq ||\nabla^2 f||_{L^2(\gamma)}$
and
\begin{equation}\label{2.29}
\lim_{\ee\ra0}||P_\ee f-f||_{\D_2^2(\gamma)}=0.
\end{equation}

Now we use $P_\ee$ to regularize $V$ and $W$. Let 
\begin{equation*}
V_n=\chi_n\,P_{1\over n}V+\log\int_{\R^d}e^{-\chi_n\,P_{1\over n}V}d\gamma\,,\quad W_n=P_{1\over n}W+\log\int_{\R^d}e^{-P_{1\over n}W}\,d\gamma,
\end{equation*}

where $\chi_n\in C_c^\infty(\R^d)$ is a smooth function with compact support satisfying usual conditions: $0\leq\chi_n\leq 1$ and
\begin{equation*}
\chi_n(x)=1\ \hbox{if}\ |x|\leq n,\quad \chi_n(x)=0\ \hbox{if}\ |x|\geq n+2,\quad \sup_{n\geq 1}||\nabla\chi_n||_\infty\leq 1.
\end{equation*}

Then the functions $V_n, W_n$ satisfy conditions in \eqref{2.28} with $2\delta_2$ for $n$ big enough,
and $\nabla V_n$ converges to $\nabla V$ in $L^2(\gamma)$. In fact,
\begin{equation*}
\nabla V_n-\nabla V=\nabla\chi_n P_{1\over n}V + \chi_n\, (\nabla P_{1\over n}V-\nabla V)+ \nabla V\,(\chi_n-1).
\end{equation*}
It is only  to check that $\dis\lim_{n\ra+\infty}\int_{\R^d} |\nabla\chi_n|^2 P_{1\over n}|V|^2\, d\gamma=0$. But 
\begin{equation*}
\int_{\R^d} |\nabla\chi_n|^2 P_{1\over n}|V|^2\, d\gamma=\int_{\R^d} |V|^2\,P_{1\over n} |\nabla\chi_n|^2\, d\gamma.\leqno(*)
\end{equation*}

For $x\in\R^d$ fixed, let $\dis r_n(x)={n-(1-e^{-1/n})|x|\over \sqrt{1-e^{-2/n}}}$, then
\begin{equation*}
P_{1\over n} |\nabla\chi_n|^2(x)\leq \int_{\R^d} {\bf 1}_{\{|e^{-1/n}x+\sqrt{1-e^{-2/n}}y|\geq n\}}d\gamma(y)
\leq \gamma(|y|\geq r_n(x))\ra 0,
\end{equation*}
as $n\ra +\infty$. Now dominated Lebesgue convergence theorem, together with above $(*)$ yields the result. 

\vskip 2mm

Let $x\ra x+\nabla\varphi_n(x)$ be the optimal transport map which pushes $e^{-V_n}\gamma$ forward to $e^{-W_n}\gamma$. By Theorem \ref{th2.4}, we have
\begin{equation}\label{2.30}
\begin{split}
&\int_{\R^d}|\nabla V_n|^2 e^{-V_n}d\gamma-\int_{\R^d}|\nabla W_n|^2 e^{-W_n}d\gamma +{2\over 1-c}\int_{\R^d}||\nabla^2W_n||_{HS}^2 e^{-W_n}d\gamma\\
&\hskip 10mm \geq 2\Ent_\gamma(e^{-V_n})-2\Ent_\gamma(e^{-W_n})+{1-c\over 2}\int_{\R^d}||\nabla^2\varphi_n||_{HS}^2 e^{-V_n}d\gamma.
\end{split}
\end{equation}
It follows that, according to \eqref{2.28},
\begin{equation*}
\sup_{n\geq 1}\int_{\R^d}||\nabla^2\varphi_n||_{HS}^2 e^{-V_n}d\gamma <+\infty.\leqno(i)
\end{equation*}

On the other hand, 
\begin{equation*}
\int_{\R^d}|\nabla\varphi_n|^2\,e^{-V_n}d\gamma=W_2^2(e^{-V_n}\gamma,e^{-W_n}\gamma).
\end{equation*}

Note that, by transport cost inequality for Guassian measure: $W_2^2(e^{-V_n}\gamma,\gamma)\leq 2\Ent_\gamma(e^{-V_n})$, the right hand side of
above equality is dominated by $ 4 (\Ent_\gamma(e^{-V_n})+ \Ent_\gamma(e^{-W_n}))$ which is bounded 
with respect to $n$, due to \eqref{2.28}. Therefore
\begin{equation*}
\sup_{n\geq 1}\int_{\R^d}|\nabla\varphi_n|^2 e^{-V_n}d\gamma <+\infty.\leqno(ii)
\end{equation*}

For the moment, we suppose that
\begin{equation*}
0<\delta_1\leq e^{-V}. \leqno(H)
\end{equation*}

Under $(H)$, above $(i)$, $(ii)$ imply that 
\begin{equation*}
\sup_{n\geq 1}\Bigl[ \int_{\R^d}|\nabla\varphi_n|^2 d\gamma +\int_{\R^d}||\nabla^2\varphi_n||_{HS}^2 d\gamma \Bigr]<+\infty.
\end{equation*}

Now by Poincar\'e inequality $\int_{\R^d}|\varphi_n-\E(\varphi_n)|^2\,d\gamma\leq\int_{\R^d}|\nabla\varphi_n|^2 d\gamma$ where
$\E(\varphi_n)$ denotes the integral of $\varphi_n$ with respect to $\gamma$. Up to changing $\varphi_n$  by $\varphi_n-\E(\varphi_n)$, we get
\begin{equation}\label{2.31}
\sup_{n\geq 1}||\varphi_n||_{\D_2^2(\gamma)}<+\infty.
\end{equation}

Therefore there exists $\varphi\in\D_2^2(\gamma)$ such that $\varphi_n\ra\varphi, \nabla\varphi_n\ra\nabla\varphi$ 
and $\nabla^2\varphi_n\ra\nabla^2\varphi$ weakly in $L^2(\gamma)$.
Now by Theorem 2.6 (for $p=1$), there exists a constant $K>0$ (independent of $n$), such that 
\begin{equation}\label{2.32}
||\nabla^2\varphi_n-\nabla^2\varphi_m||_{L^1(\gamma)}^2
\leq K\,\Bigl(||V_n-V_m||_{L^1(\gamma)}+||\nabla W_n-\nabla W_m||_{L^2(\gamma)}^2\Bigr)\ra 0,
\end{equation}
as $n,m\ra +\infty$. Also by \eqref{2.14},

\begin{equation}\label{2.33}
||\nabla\varphi_n-\nabla\varphi_m||_{L^2(\gamma)}^2
\leq  {4\over 1-c}||V_n-V_m||_{L^1(\gamma)}+ {4\over (1-c)^2}||\nabla W_n-\nabla W_m||_{L^2(\gamma)}^2\ra 0,
\end{equation}
as $n,m\ra +\infty$. It follows that $\nabla^2\varphi_n$ converges to $\nabla^2\varphi$ in $L^1(\gamma)$ and $\nabla\varphi_n$ converges
to $\nabla\varphi$ in $L^2(\gamma)$, as $n\ra +\infty$. Up to a subsequence, $\nabla^2\varphi_n$ converges to $\nabla^2\varphi$ and
$\nabla\varphi_n$ converges to $\nabla\varphi$ almost everwhere. Therefore $x+\nabla\varphi(x)$ pushes $e^{-V}\gamma$ to $e^{-W}\gamma$
and $\Id+\nabla^2\varphi$ is positive. 

\begin{theorem}\label{th2.7} Let $V\in\D_1^2(\R^d,\gamma)$ and $W\in\D_2^2(\R^d,\gamma)$ satisfying conditions \eqref{2.28} 
and $(H)$, then
the optimal transport map $x\ra x+\nabla\varphi(x)$ which pushes $e^{-V}\gamma$ to $e^{-W}\gamma$ 
is such that $\varphi\in\D_2^2(\R^d,\gamma)$ and
\begin{equation}\label{2.34}
\begin{split}
&\int_{\R^d}|\nabla V|^2 e^{-V}d\gamma-\int_{\R^d}|\nabla W|^2 e^{-W}d\gamma +{2\over 1-c}\int_{\R^d}||\nabla^2W||_{HS}^2 e^{-W}d\gamma\\
&\hskip 10mm \geq 2\Ent_\gamma(e^{-V})-2\Ent_\gamma(e^{-W})+{1-c\over 2}\int_{\R^d}||\nabla^2\varphi||_{HS}^2 e^{-V}d\gamma.
\end{split}
\end{equation}
\end{theorem}
\vskip 2mm

{\bf Proof.} Again due to \eqref{2.28}, as $n\ra+\infty$, at least for a subsequence, 
\begin{equation*}
\int_{\R^d}|\nabla V_n|^2 e^{-V_n}d\gamma\ra \int_{\R^d}|\nabla V|^2 e^{-V}d\gamma,\quad
 \int_{\R^d}|\nabla W_n|^2 e^{-W_n}d\gamma\ra \int_{\R^d}|\nabla W|^2 e^{-W}d\gamma.
 \end{equation*}

On the other hand, for  a almost everywhere convergence  subsequence, by Fatou lemma, 
$$\dis \lim_{n\ra+\infty}\int_{\R^d}||\nabla^2\varphi_n||_{HS}^2 e^{-V_n}d\gamma\geq \int_{\R^d}||\nabla^2\varphi||_{HS}^2 e^{-V}d\gamma.$$
At the limit, \eqref{2.30} leads to \eqref{2.34}.\fin

\vskip 2mm
In what follows, we will drop the condition $(H)$, but assume \eqref{2.28}.  Let $n\geq 1$, consider 
$$\dis V_n=V\wedge n.$$
Then $V_n\leq V$, $|\nabla V_n|\leq |\nabla V|$ and $V_n$ converge to $V$ in $\D_1^2(\R^d,\gamma)$.  Let $a_n=\int_{\R^d} e^{-V_n}d\gamma$; 
then $a_n\ra 1$, as $n\ra+\infty$.
Let $x\ra x+\nabla\varphi_n(x)$ be the 
optimal map which pushes $e^{-V_n}/a_n\,d\gamma$ forward to $e^{-W}d\gamma$. Then by \eqref{2.34},

\begin{equation*}
{1-c\over 2}\int_{\R^d}||\nabla^2\varphi_n||_{HS}^2 {e^{-V_n}\over a_n} d\gamma\leq
\delta_2\int_{\R^d}|\nabla V|^2\,d\gamma+{2\over 1-c}\int_{\R^d}||\nabla^2W||_{HS}^2 e^{-W}d\gamma.
\end{equation*}
 
On the other hand,
\begin{equation*}
\int_{\R^d}|\nabla\varphi_n|^2 {e^{-V}\over a_n}d\gamma\leq  \int_{\R^d}|\nabla\varphi_n|^2 {e^{-V_n}\over a_n}d\gamma=W_2^2({e^{-V_n}\over a_n} \gamma, e^{-W}\gamma).
\end{equation*}

It follows that 
\begin{equation}\label{2.35}
\sup_{n\geq 1}\Bigl[\int_{\R^d}|\nabla\varphi_n|^2 e^{-V}d\gamma+\int_{\R^d}||\nabla^2\varphi_n||_{HS}^2 e^{-V}d\gamma\Bigr]<+\infty.
\end{equation}

Since the Dirichlet form ${\cal E}(f,f)=\int_{\R^d} |\nabla f|^2\, e^{-V}d\gamma$ is closed, then there exists $Y\in \D_1^2(\R^d, \R^d; e^{-V}\gamma)$ such that
\begin{equation*}
\nabla\varphi_n\ra Y,\quad \nabla^2\varphi_n\ra \nabla Y
\end{equation*}
weakly in $L^2(e^{-V}\gamma)$.  Then, for any $\xi\in L^\infty (\R^d, \R^d; e^{-V}\gamma)$,
\begin{equation*}
\lim_{n\ra+\infty} \int_{\R^d} \langle \xi, \nabla\varphi_n\rangle\,e^{-V}d\gamma=\int_{\R^d} \langle\xi, Y\rangle\,e^{-V}d\gamma.\leqno(i)
\end{equation*}

On the other hand, by stability of optimal transport plans, there exists a $1$-convex function $\varphi\in L^1(e^{-V}\gamma)$ such that $x\ra x+\nabla\varphi(x)$
is the unique optimal transport map which pushes $e^{-V}d\gamma$ forward to $e^{-W}d\gamma$ (see \cite{Villani2},p.74), such that, up to a subsequence,
\begin{equation*}
\lim_{n\ra+\infty} \int_{\R^d} \psi(x,x+\nabla\varphi_n(x))\,{e^{-V_n}\over a_n}d\gamma
= \int_{\R^d} \psi(x,x+\nabla\varphi(x))\,e^{-V}d\gamma,\leqno(ii)
\end{equation*}
for any bounded continuous function $\psi:\R^d\times\R^d\ra\R$. Let $\alpha_R$ be a cut-off function on $\R$: $\alpha_R\in C_b(\R)$ such that $0\leq \alpha_R\leq 1$
and $\alpha_R=1$ over $[0,R]$ and $\alpha_R=0$ over $[2R,+\infty[$. Take $\xi$ as a bounded continuous function $\R^d\ra\R^d$ and consider
\begin{equation*}
\psi(x,y)=\langle\xi(x), y\rangle \alpha_R(|y|).
\end{equation*}
By above $(ii)$, and noting $\nabla\Phi_n(x)=x+\nabla\varphi_n(x)$ and $\nabla\Phi(x)=x+\nabla\varphi(x)$, we have
\begin{equation*}
\lim_{n\ra+\infty} \int_{\R^d}\langle \xi(x), \nabla\Phi_n(x)\rangle \alpha_R(|\nabla\Phi_n(x)|){e^{-V_n}\over a_n}d\gamma
=\int_{\R^d}\langle\xi(x), \nabla\Phi(x)\rangle \alpha_R(|\nabla\Phi(x)|)e^{-V}d\gamma.\leqno(iii)
\end{equation*}

Note that
\begin{equation*}
\begin{split}
& \Bigl| \int_{\R^d}\langle \xi(x), \nabla\Phi_n(x)\rangle \bigl(1-\alpha_R(|\nabla\Phi_n(x)|)\bigr)\,{e^{-V_n}\over a_n}d\gamma\Bigr|\\
&= \Bigl|\int_{\R^d} \langle \xi((\nabla\Phi_n)^{-1}(y)), y\rangle \bigl(1-\alpha_R(|y|)\bigr)\,e^{-W}d\gamma\Bigr|
\leq \delta_2\,||\xi||_\infty\,\int_{\{|y|\geq R\}}|y|\, d\gamma(y),
\end{split}
\end{equation*}

Combining this estimate with above $(iii)$, we get
\begin{equation}\label{2.36}
\lim_{n\ra+\infty} \int_{\R^d}\langle \xi(x), \nabla\Phi_n(x)\rangle \,{e^{-V_n}\over a_n}d\gamma
=\int_{\R^d}\langle\xi(x), \nabla\Phi(x)\rangle \,e^{-V}d\gamma.
\end{equation}

From \eqref{2.36}, it is not hard to see that 
\begin{equation*}
\lim_{n\ra+\infty} \int_{\R^d}\langle \xi(x), \nabla\Phi_n(x)\rangle \,e^{-V}d\gamma
=\int_{\R^d}\langle\xi(x), \nabla\Phi(x)\rangle \,e^{-V}d\gamma.
\end{equation*}
 
Now comparing with $(i)$, we get that $\nabla\Phi(x)=x+Y(x)$ or $Y=\nabla\varphi$. 

\begin{theorem}\label{th2.8} 
 Let $V\in\D_1^2(\R^d,\gamma)$ and $W\in\D_2^2(\R^d,\gamma)$ satisfying conditions \eqref{2.28}.  Then
the optimal transport map $x\ra x+\nabla\varphi(x)$ which pushes $e^{-V}\gamma$ to $e^{-W}\gamma$ 
is such that $\varphi\in\D_2^2(\R^d,\gamma)$ and
\begin{equation*}
\begin{split}
&\int_{\R^d}|\nabla V|^2 e^{-V}d\gamma-\int_{\R^d}|\nabla W|^2 e^{-W}d\gamma +{2\over 1-c}\int_{\R^d}||\nabla^2W||_{HS}^2 e^{-W}d\gamma\\
&\hskip 10mm \geq 2\Ent_\gamma(e^{-V})-2\Ent_\gamma(e^{-W})+{1-c\over 2}\int_{\R^d}||\nabla^2\varphi||_{HS}^2 e^{-V}d\gamma.
\end{split}
\end{equation*}
\end{theorem}

\vskip 2mm
{\bf Proof.} Replacing $V$ by $V_n$ in \eqref{2.34} and note that
\begin{equation*}
\underline{\lim}_{n\ra+\infty}\int_{\R^d}||\nabla^2\varphi_n||_{HS}^2\,{e^{-V_n}\over a_n}d\gamma
\geq \underline{\lim}_{n\ra+\infty}\int_{\R^d}||\nabla^2\varphi_n||_{HS}^2\,{e^{-V}\over a_n}d\gamma
\geq\int_{\R^d}||\nabla^2\varphi||_{HS}^2\,e^{-V}d\gamma,
\end{equation*}
we get the result by letting $n\ra +\infty$ in \eqref{2.34}. It remains to prove that $\varphi\in L^2(e^{-V}\gamma)$. In fact, let $\Gamma_0$ be the
optimal plan induced by $x\ra x+\nabla\varphi(x)$. Then (see section 1), under $\Gamma_0$,
\begin{equation*}
\varphi(x)+\psi(y)=|x-y|^2.
\end{equation*}

But we have seen in section 1 that $\psi\in L^2(e^{-W}\gamma)$. Then under $\Gamma_0$,
\begin{equation*}
\varphi(x)^2\leq 2\psi(y)^2 + 2 |x-y|^4.
\end{equation*}

Let $\Omega$ be the set of couples $(x,y)$ such that above inequality holds, then $\Gamma_0(\Omega)=1$. We have
\begin{equation*}
\int_{\R^d\times\R^d}\varphi^2\,d\Gamma_0=\int_{\Omega} \varphi^2\,d\Gamma_0\leq 2\int_{\R^d} \psi^2 d\Gamma_0 +2 \int_{\R^d\times\R^d}|x-y|^4\,d\Gamma_0(x,y).
\end{equation*}

It follows that
\begin{equation*}
\int_{\R^d} \varphi^2\,e^{-V}d\gamma\leq 2\int_{\R^d}\psi^2\,e^{-W}d\gamma + 16\delta_2 \int_{\R^d} |x|^4 d\gamma(x),
\end{equation*}
which is finite. The proof is complete. \fin

\vskip 2mm
We conclude this section by the following result.

\begin{theorem}\label{th2.9} Let $V_1,V_2\in\D_1^2(\R^d,\gamma)$ and $W_1,W_2\in\D_2^2(\R^d,\gamma)$ satisfying \eqref{2.28} and $(H)$.
Let $\nabla\varphi_1,  \nabla\varphi_2$ be the associated optimal transport maps. Then  for $1\leq p<2$ 
\begin{equation}\label{2.37}
\begin{split}
||\nabla^2\varphi_1-\nabla^2\varphi_2||_{L^p(e^{-V_2}\gamma)}^2 
\leq &2M(\nabla^2\varphi_1,\nabla^2\varphi_2) \Bigl[3\int_{\R^d} (V_1-V_2)e^{-V_2}d\gamma \\
&+ {2\over1- c}\int_{\R^d}|\nabla(W_1-W_2)|^2e^{-W_2}d\gamma\Bigr],
\end{split}
\end{equation}
where
  \begin{equation*}
 M(\nabla^2\varphi_1,\nabla^2\varphi_2)
 =\max\Bigl(\Bigl\|||I+\nabla^2\varphi_1||_{op}\Bigr\|_{L^{2p\over2-p}(e^{-V_2}\gamma)}^2, 
 \Bigl\|||I+\nabla^2\varphi_2||_{op}\Bigr\|_{L^{2p\over2-p}(e^{-V_2}\gamma)}^2\Bigr).
 \end{equation*}
\end{theorem}

\section{Monge-Amp\`ere equations on the Wiener space}

\quad Let's begin with finite dimension case. 

\subsection{Monge-Amp\`ere equations in finite dimension}

\begin{theorem}\label{th3.1}
 Let $V\in\D_1^2(\R^d,\gamma)$ and $W\in\D_2^2(\R^d,\gamma)$ satisfying conditions \eqref{2.28} and $(H)$. Then the
 optimal transport map $x\ra x+\nabla\varphi(x)$ from $e^{-V}\gamma$ to $e^{-W}\gamma$ solves the following 
 Monge-Amp\`ere equation
 \begin{equation}\label{3.1}
 e^{-V}=e^{-W(\nabla\Phi)}e^{\L\varphi_-{1\over2}|\nabla\varphi|^2}\det_2(\Id+\nabla^2\varphi),
 \end{equation}
 where $\nabla\Phi(x)=x+\nabla\varphi(x)$. 
 \end{theorem}
 \vskip 2mm
 
 {\bf Proof.} Let $V_n, W_n$ be the 
approximating sequence considered in section 2.2. Then  
\begin{equation}\label{3.2}
e^{-V_n}=e^{-W_n(\nabla\Phi_n)}e^{\L\varphi_n-{1\over2}|\nabla\varphi_n|^2}\det_2(\Id+\nabla^2\varphi_n),
\end{equation}
where $\nabla\Phi_n(x)=x+\nabla\varphi_n(x)$ is the optimal mal pushing $e^{-V_n}\gamma$ forward to 
$e^{-W_n}\gamma$. In order to pass to the limit in  \eqref{3.2}, we have to prove the convergence of 
$\L\varphi_n$ to $\L\varphi$, and $W_n(\nabla\Phi_n)$ to $W(\nabla\Phi)$. By \eqref{2.31}-\eqref{2.33}, we see that for any $1<p<2$,
up to a subsequence
\begin{equation*}
\lim_{n\ra +\infty}||\varphi_n-\varphi||_{\D_2^p(\gamma)}=0.
\end{equation*}

Now by Meyer inequality for Gaussian measure (see \cite{Malliavin}),
\begin{equation*}
\int_{\R^d}|\L\varphi_n-\L\varphi|^p\, d\gamma\leq C_p\, ||\varphi_n-\varphi||_{\D_2^p(\gamma)}^p.
\end{equation*}
Therefore for a subsequence, $\L\varphi_n\ra \L\varphi$ almost all. Now 

\begin{equation}\label{3.3}
\int_{\R^d} |W_n(\nabla\Phi_n)-W(\nabla\Phi)|\,d\gamma
\leq \int_{\R^d} |W_n(\nabla\Phi_n)-W(\nabla\Phi_n)|\,d\gamma+ \int_{\R^d} |W(\nabla\Phi_n)-W(\nabla\Phi)|\,d\gamma.
\end{equation}

By condition $(H)$, the first term of the right hand side of \eqref{3.3} is less than
\begin{equation*}
{1\over\delta_1}\int_{\R^d} |W_n(\nabla\Phi_n)-W(\nabla\Phi_n)|\,e^{-V_n}d\gamma
={1\over\delta_1}\int_{\R^d} |W_n-W|\,e^{-W_n}d\gamma\ra 0,
\end{equation*}
as $n\ra +\infty$. For estimating the second term, let $\ee>0$, choose $\hat W\in C_b(\R^d)$ such that 
$$||W-\hat W||_{L^1(\gamma)}\leq\ee.$$
We have

\begin{equation*}
\begin{split}
\int_{\R^d}|W(\nabla\Phi_n)-W(\nabla\Phi)|\,d\gamma
&\leq {1\over\delta_1} \int_{\R^d}|W-\hat W|(\nabla\Phi_n)\,e^{-V_n}d\gamma\\
&\hskip -15mm +\int_{\R^d}|\hat W(\nabla\Phi_n)-\hat W(\nabla\Phi)|\,d\gamma+{1\over\delta_1}\int_{\R^d}|W-\hat W|(\nabla\Phi)\,e^{-V}d\gamma\\
&\hskip -15mm \leq {2\delta_2\over\delta_1}||W-\hat W||_{L^1(\gamma)} + \int_{\R^d}|\hat W(\nabla\Phi_n)-\hat W(\nabla\Phi)|\,d\gamma.
\end{split}
\end{equation*}
It follows that 
\begin{equation*}
\lim_{n\ra+\infty}\int_{\R^d}|W(\nabla\Phi_n)-W(\nabla\Phi)|\,d\gamma=0.
\end{equation*}
So, combining this with \eqref{3.3}, up to a subsequence, $W_n(\nabla\Phi_n)\ra W(\nabla\Phi)$ almost all. The proof of \eqref{3.1} is complete. \fin

\vskip 2mm
\quad In what follows, we will drop the condition $(H)$.  

\begin{theorem}\label{th3.2}
Under conditions in Theorem \ref{th2.8}, then $\L\varphi$ exists in $L^1(\R^d,e^{-V}d\gamma)$ and 
 \begin{equation*}\
 e^{-V}=e^{-W(\nabla\Phi)}e^{\L\varphi_-{1\over2}|\nabla\varphi|^2}\det_2(\Id+\nabla^2\varphi),
 \end{equation*}
 where $\nabla\Phi(x)=x+\nabla\varphi(x)$. 
\end{theorem}

\vskip 2mm
{\bf Proof.}  Consider $V_n=V\wedge n$ for $n\geq 1$; then $V_m\leq V_n$ if $m\leq n$. Set $a_n=\int_{\R^d} e^{-V_n}\,d\gamma$, which goes to $1$ as $n\ra +\infty$.
Without loss of generality, we assume that ${1\over 2}\leq a_n\leq 2$. Let $x\ra x+\varphi_n(x)$ 
be the optimal map from ${e^{-V_n}\over a_n}d\gamma$ to $e^{-W}d\gamma$. By Theorem \ref{th2.7} or Theorem \ref{th2.8}, 
\begin{equation*}
\int_{\R^d} ||\Id +\nabla^2\varphi_n||_{op}^2\,{e^{-V_n}\over a_n}\,d\gamma
\leq 2\Bigl( 1 + {2\over 1-c}\int_{\R^d} |\nabla V_n|^2 {e^{-V_n}\over a_n}d\gamma 
+({2\over 1-c})^2\int_{\R^d} ||\nabla^2 W||_{HS}^2 e^{-W}d\gamma\Bigr),
\end{equation*}
and 

\begin{equation*}
\begin{split}
&\int_{\R^d} ||\Id +\nabla^2\varphi_m||_{op}^2\,{e^{-V_n}\over a_n}\,d\gamma
\leq 2\int_{\R^d} \Bigl( 1+ ||\nabla^2\varphi_m||_{HS}^2\Bigr)\,{e^{-V_m}\over a_m}\,e^{V_m-V_n}\,{a_m\over a_n}\,d\gamma\\
&\hskip 15mm\leq 8 \int_{\R^d} \Bigl( 1+ ||\nabla^2\varphi_m||_{HS}^2\Bigr)\,{e^{-V_m}\over a_m}\,d\gamma\\
&\hskip 15mm \leq 8\Bigl( 1 + {2\over 1-c}\int_{\R^d} |\nabla V_m|^2 {e^{-V_m}\over a_m}d\gamma 
+({2\over 1-c})^2\int_{\R^d} ||\nabla^2 W||_{HS}^2 e^{-W}d\gamma\Bigr).
\end{split}
\end{equation*}

Therefore according to Thorem \ref{th2.9}, it exists a constant $C>0$ independent of $n$, such that 
\begin{equation*}\
{1\over a_n}\int_{\R^d} ||\nabla^2\varphi_n-\nabla^2\varphi_m||_{HS}\,e^{-V}d\gamma
\leq C\, \int_{\R^d} |V_n-V_m|\,{e^{-V_n}\over a_n}d\gamma \leq 2C\delta_2 ||V_n-V_m||_{L^2(\gamma)}.
\end{equation*}

It follows that $\{\nabla^2\varphi_n;\ n\geq\}$ is a Cauchy sequence in $L^1(e^{-V}d\gamma)$. Up to subsequence, $\nabla^2\varphi_n$ converges
to $\nabla^2\varphi$ almost all. On the other hand, by Theorem \ref{th2.1},
\begin{equation*}
\int_{\R^d} |\nabla\varphi_n-\nabla\varphi_m|^2\,{e^{-V_n}\over a_n}\,d\gamma
\leq {4\over 1-c} \int_{\R^d} |V_n-V_m+\log{a_n}-\log{a_m}|\, {e^{-V_n}\over a_n}\,d\gamma,
\end{equation*}
which tends to $0$ as $m,n\ra +\infty$. Therefore up to a subsequence, $\nabla\varphi_n$ converges to $\nabla\varphi$ almost all. 

\quad Now using Theorem \ref{th3.1}, we have

\begin{equation}\label{3.4}
 {e^{-V_n}\over a_n}=e^{-W(\nabla\Phi_n)}e^{\L\varphi_n-{1\over2}|\nabla\varphi_n|^2}\det_2(\Id+\nabla^2\varphi_n),
 \end{equation}
 where $\nabla\Phi_n(x)=x+\nabla\varphi_n(x)$.  As what did in the last part of the proof to Theorem \ref{th3.1}, we have
 \begin{equation}\label{3.5}
 \lim_{n\ra \infty} \int_{\R^d} |e^{-W(\nabla\Phi_n)}-e^{-W(\nabla\Phi)}|\, e^{-V}d\gamma=0.
 \end{equation}
Therefore for a subsequence, we proved that each term except $\L\varphi_n$ in \eqref{3.4} converges almost all; it follows 
\begin{equation}\label{3.6}
\hbox{up to a subsequence,} \ \L\varphi_n\ \hbox{converges to a function } F\ \hbox{almost all}.
\end{equation}
The fact that $F\in L^1(\R^d, e^{-V}d\gamma)$ comes from the relation

\begin{equation*}
F= -V + W(\nabla\Phi)+{1\over 2}|\nabla\varphi|^2-\log{\det_2(\Id+\nabla^2\varphi)}.
\end{equation*}

Now it remains to prove that $\L\varphi$ exists in $L^1(\R^d,e^{-V}d\gamma)$ and $F=\L\varphi$. 
The difficulty is that we have no more the control in $L^2(e^{-V}d\gamma)$ of $\L\varphi_n$
by $\nabla^2\varphi_n$. We will proceed as in \cite{BogachevKolesnikov}.

\begin{lemma}\label{th3.3}
Assume that $e^{-V}\geq\delta_1>0$. Then there exists a constant $K$ independent of $\delta_1$ such that
for any $f\in \D_2^2(\R^d, e^{-V}d\gamma)$,
\begin{equation}\label{3.7}
\int_{\R^d}(\L f)^2 e^{-|\nabla f|^2}\,e^{-V}d\gamma
\leq K\,\Bigl(1+\int_{\R^d}|\nabla^2 f|^2\, e^{-V}d\gamma +\int_{\R^d}|\nabla V|^2\,e^{-V}d\gamma\Bigr).
\end{equation}
\end{lemma}

\vskip 2mm

{\bf Proof.} Any $f\in \D_2^2(\R^d, e^{-V}d\gamma)$ is also in $\D_2^2(\R^d, d\gamma)$; then $\L f$ exists in $L^2(\R^d, e^{-V}d\gamma)$, and we can
approximate $f$ by functions in $C^2$ bounded with bounded derivatives up to order $2$. For the moment, assume that $f$ is in the latter class. So

\begin{equation}\label{3.8}
\int_{\R^d} (\L f)^2 e^{-|\nabla f|^2}\,e^{-V}d\gamma=- \int_{\R^d} \langle\nabla f,\, \nabla(\L f e^{-|\nabla f|^2}e^{-V})\rangle\,d\gamma.
\end{equation}
We have
\begin{equation}\label{3.9}
\begin{split}
&\langle\nabla f,\, \nabla(\L f e^{-|\nabla f|^2}e^{-V})\rangle
=\langle\nabla f, \nabla\L f\rangle\, e^{-|\nabla f|^2}e^{-V}\\
&\hskip 10mm -2\langle \nabla f\otimes\nabla f, \nabla^2 f\rangle e^{-V}\L f e^{-|\nabla f|^2}
-\langle \nabla f,\nabla V\rangle\L f e^{-|\nabla f|^2} e^{-V}.
\end{split}
\end{equation}

By Cauchy-Schwarz inequality,

\begin{equation*}
\begin{split}
&\int_{\R^d} \langle \nabla f\otimes\nabla f, \nabla^2 f\rangle e^{-V}\L f e^{-|\nabla f|^2}d\gamma\\
&\hskip -15mm \leq \Bigl(\int_{\R^d} \langle \nabla f\otimes\nabla f, \nabla^2 f\rangle^2 e^{-|\nabla f|^2}e^{-V}d\gamma\Bigr)^{1/2}
\Bigl(\int_{\R^d} (\L f)^2 e^{-|\nabla f|^2}e^{-V}d\gamma\Bigr)^{1/2}.
\end{split}
\end{equation*}

In the same way, we treat the last term in \eqref{3.9}. Set $A= \int_{\R^d} \langle\nabla f,\,\nabla\L f\rangle e^{-|\nabla f|^2}e^{-V}d\gamma$, 
\begin{equation*}
B= 2 \Bigl(\int_{\R^d}  \langle \nabla f\otimes\nabla f, \nabla^2 f\rangle^2 e^{-|\nabla f|^2}e^{-V}d\gamma\Bigr)^{1/2}
+ \Bigl(\int_{\R^d} \langle \nabla f,\nabla V\rangle^2 e^{-|\nabla f|^2} e^{-V}d\gamma\Bigr)^{1/2},
\end{equation*}
and $Y=\Bigl(\int_{\R^d} (\L f)^2 e^{-|\nabla f|^2}e^{-V}d\gamma\Bigr)^{1/2}$. Then combining \eqref{3.8}, \eqref{3.9} and par above computation, we get

\begin{equation}\label{3.10}
Y^2\leq -A + BY.
\end{equation}

It follows that the discriminant of $P(\lambda)=\lambda^2 - B\lambda + A$ is non negative 
and $P(\lambda)=(\lambda-\lambda_1)(\lambda-\lambda_2)$. The relation \eqref{3.10} implies that $Y$ is between two roots of $P$. In particular,

\begin{equation}\label{3.11}
Y\leq (B+\sqrt{B^2 - 4A})/2.
\end{equation}

It is obvious that for a numerical constant $K_1>0$, 
\begin{equation*}
B^2 \leq K_1\,\Bigl(\int_{\R^d}|\nabla^2 f|^2\, e^{-V}d\gamma +\int_{\R^d}|\nabla V|^2\,e^{-V}d\gamma\Bigr).
\end{equation*}

For estimating the term $A$, we use the commutation formula for Gaussian measures (see for example \cite{Fang}, p. 144),
\begin{equation*}
\nabla\L f = \L\nabla f -\nabla f,
\end{equation*}

so that we get 
\begin{equation*}
|A|\leq K_1\,\Bigl(1+\int_{\R^d}|\nabla^2 f|^2\, e^{-V}d\gamma +\int_{\R^d}|\nabla V|^2\,e^{-V}d\gamma\Bigr).
\end{equation*}

Now the relation \eqref{3.11} yields \eqref{3.7}. \fin

\vskip 2mm
Applying \eqref{3.7} to $\varphi_n$, we have
\begin{equation*}
\sup_{n\geq 1} \int_{\R^d} (\L\varphi_n)^2 e^{-|\nabla\varphi_n|^2}\, e^{-V}d\gamma <+\infty.
\end{equation*}

Therefore the family $\{\L\varphi_n\, e^{-|\nabla\varphi_n|^2/2}\}$ is uniformly integrable with respect to $e^{-V}d\gamma$. Then for any $\xi\in C_b^1(\R^d)$,

\begin{equation}\label{3.12}
\lim_{n\ra +\infty} \int_{\R^d} \L\varphi_n e^{-|\nabla\varphi_n|^2/2}\xi\,e^{-V}d\gamma
= \int_{\R^d} F e^{-|\nabla\varphi|^2/2}\xi\,e^{-V}d\gamma.
\end{equation}

But 
\begin{equation*}
\begin{split}
\int_{\R^d} \L\varphi_n e^{-|\nabla\varphi_n|^2/2}\xi\,e^{-V}d\gamma
=&\int_{\R^d}\langle\nabla\varphi_n\otimes \nabla\varphi_n,\,\nabla^2\varphi_n\rangle e^{-|\nabla\varphi_n|^2/2}\,\xi e^{-V}d\gamma\\
&-\int_{\R^d} \langle\varphi_n, \nabla(\xi e^{-V})\rangle e^{-|\nabla\varphi_n|^2/2}d\gamma,
\end{split}
\end{equation*}

which converges to $\int_{\R^d}\langle\nabla\varphi\otimes \nabla\varphi,\,\nabla^2\varphi\rangle e^{-|\nabla\varphi|^2/2}\,\xi e^{-V}d\gamma
-\int_{\R^d} \langle\varphi, \nabla(\xi e^{-V})\rangle e^{-|\nabla\varphi|^2/2}d\gamma$. So we get

\begin{equation}\label{3.13}
\int_{\R^d} (F-\langle\nabla\varphi,\nabla V\rangle) e^{-|\nabla\varphi|^2/2}\xi\,e^{-V}d\gamma
=-\int_{\R^d} \langle \nabla\varphi,\, \nabla(\xi e^{-|\nabla\varphi|^2/2})\rangle\, e^{-V}d\gamma.
\end{equation}

Note that the generator $\L_V$ associated to the Dirichlet form ${\cal E}_V(f,f)=\int_{\R^d} |\nabla f|^2\, e^{-V}d\gamma$ admits the expression 
$\L_V(f)=\L(f)-\langle \nabla f,\nabla V\rangle$. Therefore the relation \eqref{3.13} tells us that $F=\L\varphi$. \fin

\subsection{Monge-Amp\`ere equations on the Wiener space}

 \quad We return now to the situation in Theorem \ref{th1.1}. Let $V\in\D_1^2(X)$ and  $W\in\D_2^2(X)$ such that $\int_X e^{-V}d\mu=\int_X e^{-W}d\mu=1$.
Assume that
\begin{equation}\label{3.14}
e^{-V}\leq \delta_2,\quad e^{-W}\leq \delta_2,
\end{equation}
and \eqref{1.8}. Let $\{e_n; n\geq 1\}\subset X^*$ be an orthonormal basis of $H$ and $H_n$ the subspace spanned
by $\{e_1, \ldots, e_n\}$. As in section 1, denote  $\dis\pi_n(x)=\sum_{j=1}^n e_j(x)e_j$ and $\F_n$ the sub $\sigma$-field generated
by $\pi_n$. 
In the sequel, we will see that the manner to regularize the density functions $e^{-V}$ and $e^{-W}$ has impacts on final results.
\vskip 2mm

 Set
\begin{equation}\label{3.15}
\E(e^{-V}|\F_n)=e^{-V_n}\circ\pi_n,\ \E(W|\F_n)=W_n\circ\pi_n.
\end{equation}

It is obvious that $\nabla^2 W_n\geq -c\, \Id_{H_n\otimes H_n}$. Applying Theorem \ref{th2.8}, there is a $\varphi_n\in \D_2^2(H_n,\gamma_n)$
such that $x\ra x+\nabla\varphi_n(x)$ is the optimal transport map
 which pushes $e^{-V_n}\gamma_n$ to $e^{-W_n}\gamma_n$. Let $\tilde \varphi_n=\varphi_n\circ\pi_n$.
 We have
 
 \begin{equation}\label{3.16}
 \begin{split}
 &{1-c\over 2}\int_{H_n} ||\nabla^2\varphi_n||_{HS}^2 e^{-V_n}d\gamma_n\\
 &\hskip -5mm\leq \int_{H_n}|\nabla V_n|^2 e^{-V_n}d\gamma_n + {2\over 1-c}\int_{H_n}||\nabla^2 W_n||_{HS}^2 e^{-W_n}d\gamma_n.
 \end{split}
 \end{equation}
 
By Cauchy-Schwarz inequality for conditional expectation, 
\begin{equation*}
|\nabla \E(e^{-V}|\F_n)|_{H_n}^2\leq \E(|\nabla V|_H^2 e^{-V}|\F_n)\,\E(e^{-V}|\F_n)
\end{equation*}
which implies that $\int_{H_n} |\nabla V_n|^2 e^{-V_n}d\gamma_n \leq \int_X |\nabla V|^2 e^{-V}d\mu$. So \eqref{3.16} yields 

 \begin{equation}\label{3.17}
{1-c\over 2}\int_{X} ||\nabla^2\tilde\varphi_n||_{HS}^2 e^{-V}d\mu
 \leq \int_{X}|\nabla V|^2 e^{-V}d\mu + {2\delta_2\over 1-c}\int_{X}||\nabla^2 W||_{HS}^2 d\mu.
 \end{equation}

Let $n,m$ be two integers such that $n>m$, and $\pi_m^n: H_n\ra H_m$ the orthogonal projection.   Then
$I_{H_n}+\nabla(\varphi_m\circ\pi_m^n)$ pushes $e^{-V_m}\circ\pi_m^n\gamma_n$ to $e^{-W_m}\circ\pi_m^n\,\gamma_n$. In fact, 
for any bounded continuous function $f: H_n\ra\R$, 
\begin{equation*}
\int_{H_n} f \bigl(x+\pi_m^n(\nabla\varphi_m)\circ\pi_m^n(x)\bigr)e^{-V_m}\circ\pi_m^nd\gamma_n
=\int_{H_m^\perp}\Bigl[\int_{H_m} f(z'+z+\pi_m^n(\nabla\varphi_m)(z))e^{-V_m}(z)d\gamma_m(z)\Bigr]d\hat\gamma(z'),
\end{equation*}
where $H_n=H_m\oplus H_m^\perp$ and $\gamma_n=\gamma_m\otimes\hat\gamma$. 
Note that $\pi_m^n(\nabla\varphi_m)=\nabla\varphi_m$; then the last term in above equality yields
\begin{equation*}
\int_{H_m^\perp}\Bigl[\int_{H_m} f(z'+y)e^{-W_m}(y)d\gamma_m(y)\Bigr]d\hat\gamma(z')
=\int_{H_n} f(x)e^{-W_m}\circ \pi_m^n(x)d\gamma_n(x).
\end{equation*}

Now by \eqref{2.14},
\begin{equation*}
\begin{split}
&||\nabla\varphi_n-\nabla(\varphi_m\circ\pi_m^n)||_{L^2(e^{-V_n}\gamma_n)}^2\\
&\hskip -10mm\leq {4\over 1-c}\int (V_n-V_m\circ\pi_m^n)e^{-V_n}d\gamma_n+{4\over (1-c)^2}\int_{H_n}|\nabla W_n-\nabla(W_m\circ\pi_m^n)|^2\, e^{-W_n}\,d\gamma_n,
\end{split}
\end{equation*}

or 
\begin{equation}\label{3.18}
\begin{split}
&||\nabla\tilde\varphi_n-\nabla\tilde\varphi_m||_{L^2(e^{-V}\mu)}^2\\
&\hskip -10mm\leq {4\over 1-c}\int_X (V_n\circ\pi_n-V_m\circ\pi_m)e^{-V}d\mu
+{4\delta_2\over (1-c)^2}\int_{X}|\nabla \E( W|\F_n)-\nabla\E(W|\F_m)|^2\,d\mu.
\end{split}
\end{equation}

Now in order to control the sequence of functions $\tilde\varphi_n$, we suppose that
\begin{equation}\label{3.19}
e^{-V}\geq\delta_1>0.
\end{equation}

Under \eqref{3.19}, it is clear that
\begin{equation*}
\int_X (V_n\circ\pi_n-V_m\circ\pi_m)e^{-V}d\mu\ra 0,\ \hbox{as }n,m\ra +\infty.
\end{equation*}

Now replacing $\tilde\varphi_n$ by $\tilde\varphi_n-\int_X \tilde\varphi_n\,d\mu$ and according to Poincar\'e inequality,
and by \eqref{3.18}, we see that $\tilde\varphi_n$ converges in $\D_1^2(X)$ to a function $\varphi$. On the other hand, 
by \eqref{3.17}, $\tilde\varphi_n$ converges to a function $\hat\varphi\in \D_2^2(X)$ weakly. By uniqueness of limits, we see in fact that
$\varphi\in\D_2^2(X)$. Now we proceed as in subsection 3.1, we have
\begin{equation}\label{3.20}
\lim_{n\ra +\infty}\int_X ||\nabla^2\tilde\varphi_n-\nabla^2\varphi||_{HS}\, d\mu=0.
\end{equation}

Combining \eqref{3.20} and \eqref{3.17}, up to a subsequence, for any $1<p<2$,
\begin{equation}\label{3.21}
\lim_{n\ra +\infty}\int_X ||\nabla^2\tilde\varphi_n-\nabla^2\varphi||_{HS}^p\, d\mu=0.
\end{equation}

By Meyer inequality (\cite{Malliavin}),
\begin{equation}\label{3.22}
\lim_{n\ra +\infty}\int_X ||\L\tilde\varphi_n-\L\varphi||_{HS}^p\, d\mu=0.
\end{equation}

So everything goes well under the supplementary condition \eqref{3.19}. We finally get

\begin{theorem}\label{th3.4}
Under conditions \eqref{3.14}, \eqref{1.8} and \eqref{3.19}, there exists a function $\varphi\in\D_2^2(X)$ such that
$x\ra x+\nabla\varphi(x)$ pushes $e^{-V}\mu$ to $e^{-W}\mu$ and solves the Monge-Amp\`ere equation
\begin{equation*}
e^{-V}=e^{-W(T)}e^{\L\varphi-{1\over 2}|\nabla\varphi|^2}\,det_2(\Id_{H\otimes H}+\nabla^2\varphi),
\end{equation*}
where $T(x)=x+\nabla\varphi(x)$.
\end{theorem}

{\bf Remark:} The regularization of $W$ used in \eqref{3.15} does not allows to prove that 
$W_2^2(e^{-V_n}\gamma_n, e^{-W_n}\gamma_n)$ converges to $W_2^2(e^{-V}\mu,e^{-W}\mu)$ contrary to section 1; 
we do not know if the map $T$ constructed in Theorem \ref{th3.4} is the optimal transport : which is due to the singularity of
the cost function $d_H$ in contrast to finite dimensional case (see subsection 3.1). 

\begin{theorem}\label{th3.5}
Assume all conditions in Theorem \ref{th3.4} and that 
\begin{equation}\label{3.23}
W_n\ \hbox{defined in \eqref{1.11} is in } \D_2^2(H_n)\ \hbox{ for all } n\geq 1. 
\end{equation}
Then there is 
a function $\varphi\in\D_2^2(X)$ such that $x\ra T(x)= x+\nabla\varphi(x)$ is the optimal transport map which pushes $e^{-V}\mu$ to $e^{-W}\mu$
and $T$ is the inverse map of $S$ in Theorem \ref{th1.1}.
\end{theorem}

\vskip 2mm
{\bf Proof.} By Proposition 5.1 in \cite{FeyelUstunel3}, $W_n$ satisfies the condition \eqref{2.28}. So we can repeat the arguments as above, but the difference
is that in actual case, $W_2^2(e^{-V_n}\gamma_n, e^{-W_n}\gamma_n)$ converges to $W_2^2(e^{-V}\mu,e^{-W}\mu)$. Using notations in the proof
of Theorem \ref{th1.1}, $x\ra x-{1\over 2}\nabla\varphi_n(x)$ is the optimal transport map, which pushes $e^{-V_n}\gamma_n$ to $e^{-W_n}\gamma_n$. So that
\begin{equation*}
W_2^2(e^{-V}\mu,e^{-W}\mu)={1\over 4}\int_X|\nabla\varphi|_H^2\, e^{-V}d\mu,
\end{equation*}

that means that $x\ra T(x)= x-{1\over 2}\nabla\varphi(x)$ is the optimal transport map which pushes $e^{-V}\mu$ to $e^{-W}\mu$. To see that
$T$ is the inverse map of $S$ in Theorem \ref{1.1}, we use \eqref{1.20}, which implies that under the optimal plan $\Gamma_0$,
\begin{equation*}
-2\psi(x)+ \varphi(y)=d_H(x,y)^2,
\end{equation*}
since we have replaced $-{1\over2}\psi$ by $\psi$ at the end of the proof of Theorem \ref{th1.1}.
But now $\varphi\in\D_2^2(X)$, we can differentiate $\varphi$ as in section 1, so that under $\Gamma_0$,
\begin{equation*}
x=y-{1\over 2}\nabla\varphi(y).
\end{equation*}

Therefore $\eta\in L^2(X,H,e^{-V}\mu)$ is given by $\eta=-{1\over2}\nabla\varphi$ with $\varphi\in \D_2^2(X)$. \fin

\vskip 2mm
{\bf Examples:} (i) If $W\in \D_2^2(X)$ satisfies $\int_X |\nabla W|^4\, d\mu<+\infty$ and $0<\delta_1\leq e^{-W}\leq \delta_2$ then
\eqref{3.23} holds.\fin
\vskip 2mm

(ii) For an orthonormal basis $\{e_n;n\geq 1\}$ of $H$, define $W(x)=\sum_{n\geq 1}\lambda_n e_n(x)^2$, where $\lambda_n>-1/2$ and
$ \sum_{n\geq 1}|\lambda_n|<+\infty$. We have,
\begin{equation*}
\E(e^{-W}|\F_n)=e^{-\sum_{k=1}^n \lambda_k e_k(x)^2}\,\prod_{k>n}\E(e^{-\lambda_k e_k(x)^2})
=\alpha_n e^{-\sum_{k=1}^n \lambda_k e_k(x)^2},
\end{equation*}
where $\alpha_n=\prod_{k>n}{1\over \sqrt{1+2\lambda_k}}$. So \eqref{3.23} holds. \fin

\section{Appendix: } 

For the sake of reader's convenience, we collect in this section some results used in this work.  

\begin{lemma}\label{th4.1}
Let $A$ be a symmetric positive definite matrix and $B$ be a symmetric matrix on $\R^d$; then
\begin{equation}\label{2.22}
||A^{-1/2}BA^{-1/2}||_{HS}\geq {||B||_{HS}\over||A||_{op}},
\end{equation}
where $||\cdot||_{op}$ denotes the norm of matrices.
\end{lemma}

\vskip 2mm
{\bf Proof.} Let $C=A^{-1/2}BA^{-1/2}$, then $C=A^{1/2}BA^{1/2}$. Let $\{e_1, \cdots, e_d\}$ be an orthonormal basis of $\R^d$, of eigenvalues of $A$:
$A^{1/2}e_i=\sqrt{\lambda_i}\,e_i$.  We have $Be_i=\sqrt{\lambda_i}\,A^{1/2}Ce_i$ and
\begin{equation*}
|Be_i|^2\leq \max(\lambda_i)\,|A^{1/2}Ce_i|^2= \max(\lambda_i)\, \langle Ce_i, ACe_i\rangle\leq ||A||_{op}^2\,|Ce_i|^2.
\end{equation*}
It follows that $||B||_{HS}^2\leq ||A||_{op}^2\,||C||_{HS}^2$. The result \eqref{2.22} follows. \fin

\begin{lemma}\label{th4.2} Let $A, B$ be symmetric matrices such that $I+A$ and $I+B$ are positive definite. Then
\begin{equation}\label{2.23}
\begin{split}
&-\log\det_2\Bigl((I+A)(I+B)^{-1}\Bigr)\\
&\hskip -10mm=\int_0^1 (1-t)||(I+(1-t)B+tA)^{-1/2}(A-B)(I+(1-t)B+tA)^{-1/2}||_{HS}^2\,dt.
\end{split}
\end{equation}
\end{lemma}

\vskip 2mm

{\bf Proof.} Note first $I-(I+A)(I+B)^{-1}=(B-A)(I+B)^{-1}$ and 
\begin{equation*}
\hbox{Trace}\Bigl[I-(I+A)(I+B)^{-1}\Bigr]=\langle B-A, (I+B)^{-1}\rangle_{HS}.\leqno(i)
\end{equation*}

Let $\chi(t)=\log\det\Bigl(I+(1-t)B+tA\Bigr)$ for $t\in [0,1]$. We have
\begin{equation*}
\chi'(t)=\hbox{Trace}\Bigl[(A-B)(I+(1-t)B+tA)^{-1}\Bigr]
=\langle A-B,(I+(1-t)B+tA)^{-1}\rangle_{HS}.
\end{equation*}
Then
\begin{equation*}
\log\det(I+A)-\log\det(I+B)
=\langle A-B,\int_0^1(I+(1-t)B+tA)^{-1}\,dt\rangle_{HS}.
\end{equation*}

According to above $(i)$ and definition of $\det_2$, we get
\begin{equation*}
\begin{split}
&-\log\det_2\Bigl((I+A)(I+B)^{-1}\Bigr)
=\langle A-B, \int_0^1\Bigl[(I+B)^{-1}-(I+(1-t)B+tA)^{-1}\Bigr]\,dt\rangle_{HS}\\
&\hskip 5mm=\langle A-B, \int_0^1\Bigl[\int_0^t (I+(1-s)B+sA)^{-1}\,(A-B)\,(I+(1-s)B+sA)^{-1}ds\Bigr]dt\rangle_{HS}
\end{split}
\end{equation*}
which is equal to $\int_0^1(1-t)\langle A-B, (I+(1-t)B+tA)^{-1}\,(A-B)\,(I+(1-t)B+tA)^{-1}\rangle_{HS}\,dt$, implying \eqref{2.23}.\fin

\end{document}